%% file: main.tex
\pdfoutput=1

\documentclass[preprint]{article}

\input{preamble}

\input{newcmd}

\title{Numerical Integration over the Unit Sphere by using \std{}\footnote{This work was supported by National Natural Science Foundation of China [grant number 11301222].}}

\author
{Congpei An$^{1}$,
Siyong Chen$^{1}$\\
\\
\normalsize{$^{1}$Department of Mathematics, Jinan Universtiy,}\\
\normalsize{ Guangzhou 510632, China}\\
\\
\normalsize{E-mail: tancpei@jnu.edu.cn, chensiy98@126.com}
}

\date{}

\begin{document}


\maketitle


\input{chapter/abstract}

\input{chapter/1-introduction}

\input{chapter/2-prequisites}

\input{chapter/3-wce}

\input{chapter/4-tran}

\input{chapter/5-numres}

\input{chapter/Final-Remark}

\input{chapter/acknowledgment}


\bibliography{refer}


\end{document}

%% file: preamble.tex
\usepackage{amsfonts}
\usepackage{amsmath}
\allowdisplaybreaks[4]  
\usepackage{amsthm}
\usepackage{amssymb}

\usepackage{mathrsfs}
\usepackage{bm}
\usepackage{extarrows}
\usepackage{graphicx}

\usepackage{epstopdf}  
\usepackage{varwidth}
\usepackage{subfigure}
\usepackage[subfigure]{ccaption}
\usepackage{enumitem}
\setlist[enumerate]{leftmargin=3em}

\graphicspath{{figure/}}

\usepackage[pdfencoding=auto]{hyperref}
\usepackage{bookmark}

\usepackage[numbers, sort,comma]{natbib}

%% file: newcmd.tex
\newcommand{\dpsty}{ \displaystyle }

\newcommand{\Stwo}{ \mathbb{S}^2 }
\newcommand{\XN}{ \mathcal{X}_N }
\newcommand{\mbfxl}[1]{ \mathbf{x}_{#1} }
\newcommand{\mbfx}{ \mathbf{x} }

\newcommand{\Pt}{ \mathbb{P}_t }
\newcommand{\df}{ d }
\newcommand{\domegax}{ \,\df \omega(\mathbf{x}) }
\newcommand{\std}{spherical $t$-design}
\newcommand{\stds}{spherical $t$-designs}
\newcommand{\Std}{Spherical $t$-designs}

\newcommand{\Ct}{ \bm{C}_t }
\newcommand{\wce}{\mathrm{wce}}

\newcommand{\dist}{\mathrm{dist}}

\newcommand{\QN}{ Q_{N} }

\newcommand{\Yt}{ \bm{Y}_t }
\newcommand{\Ylk}{ Y_{\ell,k} } 

\newcommand{\inprod}[2]{ \left\langle #1 ,\, #2 \right\rangle }
\newcommand{\Tone}{\mathcal{T}_{1}}
\newcommand{\Ttwo}{\mathcal{T}_{2}}

\newcommand{\seq}[3]{ {#1}_{#2},\,\ldots ,\, {#1}_{#3} }
\newcommand{\sequ}[2]{ {#1},\,\ldots ,\, {#2} }

\newcommand{\Sob}[1]{ \mathbb{H}^{#1} }
\newcommand{\Hs}{ \Sob{s} }
\newcommand{\Ltwo}{ \mathbb{L}_{2} }

\newcommand{\Real}{\mathbb{R}}

\newcommand{\norm}[2]{ \left\Vert#1\right\Vert_{#2} }
\newcommand{\abs}[1]{ \left\vert #1 \right\vert }

\newcommand{\vectornorm}[1]{\norm{#1}{2}}

\newcommand{\eqtref}[1]{ (\ref{#1}) }
\newcommand{\figref}[1]{ Figure~\ref{#1} }
\newcommand{\tabref}[1]{ Table~\ref{#1} }
\newcommand{\secref}[1]{Section~\ref{#1} }

\newtheorem{conjecture}{Conjecture}

\newtheorem{definition}{Definition}

\newcounter{lastnote}

%% file: chapter/abstract.tex
\abstract{
  This paper studies numerical integration over the unit sphere $ \Stwo \subset \Real^{3} $ by using \std{}, which is an equal positive weights quadrature rule with polynomial precision $t$.
  We investigate two kinds of \stds{} with $t$ up to 160. One is well conditioned \stds{}(WSTD), which was proposed by \cite{an2010well} with $ N=(t+1)^{2} $.
  The other is efficient \std{}(ESTD), given by Womersley \cite{womersley2009spherical}, which is made of roughly of half cardinality of WSTD.
  Consequently, a series of persuasive numerical evidences indicates that WSTD is better than ESTD in the sense of worst-case error in Sobolev space $ \Hs(\Stwo) $.
  Furthermore, WSTD is employed to approximate integrals of various of functions, especially including integrand has a point singularity over the unit sphere and a given ellipsoid.
  In particular, to deal with singularity of integrand, Atkinson's transformation \cite{atkinson2004quadrature} and Sidi's transformation \cite{sidi2005application} are implemented with the choices of `grading parameters' to obtain new integrand which is much smoother.
  Finally, the paper presents numerical results on uniform errors for approximating representive integrals over sphere with three quadrature rules: Bivariate trapezoidal rule, Equal area points and WSTD.
}


%% file: chapter/1-introduction.tex
\section{Introduction}
We consider numerical integration over the unit sphere
  \begin{equation*}
    \Stwo := \left\{ \mathbf{x}=[x,\, y,\, z]^{T} \in \Real^{3} \,\vert\, x^{2}+y^{2}+z^{2}=1 \right\} \subset \Real^{3}.
  \end{equation*}
The exact integral of an integrable function $f$ defined on $\Stwo$ is
  \begin{equation}
    I(f) := \int_{\Stwo} f(\mathbf{x}) \domegax,
  \end{equation}
where $\domegax$ denotes the surface measure on $\Stwo$. The aim of this paper is approximate $I(f)$ by positive weight quadrature rules of the form
  \begin{equation}\label{eqt:quadrarule}
    Q_{N}(f) := \sum_{j=1}^{N} \omega_{j} f(\mathbf{x}_{j}),
  \end{equation}
where $ \mathbf{x}_{j} \in \Stwo$, $ 0 < \omega_{j} \in \Real,\; j=1,\,\ldots,\,N $. \par

As shown in literature \cite{an2010well, brauchart2014qmc, cui1997equidistribution, grabner1993spherical, sloan2004extremal} on numerical integration over sphere, there are many point sets can be used as quadrature rules.
It is merit to consider Quasi Monte Carlo(QMC) rules, or equal weight numerical integration, for functions in a Sobolev space $\Hs(\Stwo)$, with smoothness parameter $s>1$.
In particular, \cite{brauchart2014qmc} provides an emulation of \stds{}, named sequence of QMC designs.  \par

As is known, for numerical integration, sequences of \std{} enjoy the property that it convergent very fast in Sobolev spaces \cite{hesse2005worst}. 
In this paper we focus our interest on approximating integrals with the aid of \stds{}.
The concept of \std{} was introduced by Delsarte, Goethals and Seidel in 1977 \cite{delsarte1977spherical}, as following:
  \begin{definition}
    A point set $\XN=\{ \seq{\mathbf{x}}{1}{N} \} \subset \Stwo$ is called a \std{} if an equal weight quadrature rule with node set $\XN$ is exact for all polynomials $p$ with degree no more than $t$, i.e.
      \begin{equation} \label{eq:stddef}
        \frac{1}{N}\sum_{j=1}^N p(\mbfxl{j}) = \frac{1}{4\pi} \int_{\Stwo} p(\mbfx{}) \domegax \quad\quad \forall p \in \Pt,
      \end{equation}
    where $\Pt$ is the linear space of restrictions of polynomials of degree at most $t$ in $3$ variables to $\Stwo$.
  \end{definition}
\noindent In past decades, \std{} has been studied extensively 
\cite{an2010well, bannai2009survey, bondarenko2013optimal, chen2011computational, chen2006existence, sloan2009variational}.
The existence of \std{} for all $t$ only for sufficiently large $N$ was shown in \cite{seymour1984averaging}.
However, when $t$ is given, the smallest number of $N$ to construct a \std{} is still to be fixed.
Bondarenko, Radchenko and Viazovska \cite{bondarenko2013optimal} claimed that there always exist a \std{} for $N \geq ct^{2}$, but $c$ is an unknown constant.
In practice, one might have to construct \stds{} by assist of numerical computation, when $t$ is large.
To the best of our knowledge, there are many numerical methods for constructing \stds{}.
For detail, we refer \cite{an2016numerical, chen2011computational, sloan2004extremal, womersley2009spherical}. 
However, it is not easy to overcome rounding errors in computation.
Therefore, reliable numerical \std{} is cherished in numerical construction.
\par

\cite{chen2011computational} and \cite{chen2006existence} verified \std{} exist in a small neighbourhood of extremal system  on $\Stwo$.
It is worth noting that \emph{well conditioned \stds{}(WSTD)} are not only have good geometry but also good for numerical integration with $N=(t+1)^2$, that is the dimension of the linear space $\Pt$.
In \cite{an2010well}, WSTD are constructed just up to 60.
In present paper we can use WSTD for $t$ up to 160 with $N=(t+1)^2$, from a very recently work \cite{an2016numerical}.
Moreover, Womersley introduced \emph{efficient \std{}(ESTD)}, which are consist of $N \approx \frac{1}{2}t^{2} $ points \cite{womersley2009spherical}.
Both WSTD and ESTD are applied to approximate the integral of a well known smooth function -- Franke function.
High attenuation to absolute error \cite{an2010well, brauchart2014qmc, womersley2009spherical}, excellent performance enhance the attractiveness of \stds{}.
Inspired by \cite{brauchart2014qmc, hesse2005worst}, it is natural to compare the worst-case error for these two \stds{} for $t$ up to $160$.
Consequently, we will employ the lower worst-case error \std{} (actually it is WSTD) to approximate the integral of various of functions: smooth function, $C^{0}$ function, near-singular function, singular function over $\Stwo$ and ellipsoid. \par
In sequel, we provide necessary background and terminology for spherical polynomial, \std{}. 
\secref{sec:wce} introduces the concept of worst-case errors of positive weight quadratures rules on $\Stwo$. 
Immediately, a series of numerical experiments indicates that WSTD has lower worst-case error than ESTD.
Consequently, we use WSTD to evaluate the numerical integration of several kinds of test functions in below.
\secref{sec:tran} focus on how to deal with point singularity of integrands, we apply the variable transformations raised by Atkinson \cite{atkinson2004quadrature} and Sidi \cite{sidi2005application} to obtain new smoother integrands, respectively.
In \secref{sec:num_res}, we investigate three quadrature nodes: Bivariate trapezoidal rule, Equal area points and WSTD. The geometry of these quadrature nodes is compared. We also demonstrate numerical results on uniform errors for approximating integrals of a set of test functions, by using these three quadrature nodes.

%% file: chapter/2-prequisites.tex
\section{Background}\label{sec:bg}
Let $ \Ltwo := \Ltwo(\Stwo) $ be the space of square-integrable measurable functions over $\Stwo$. The Hilbert space $\Ltwo$ is endow with the inner product
  \begin{equation*}
    \inprod{f}{g}_{\Ltwo} := \int_{\Stwo} f(\mathbf{x}) \, g(\mathbf{x}) \domegax, \quad f, g \in \Ltwo .
  \end{equation*}
And the induced norm is
  \begin{equation*}
    \norm{f}{\Ltwo} := \left(\int_{\Stwo} \abs{ f(\mathbf{x}) }^{2}  \domegax \right)^{\frac{1}{2}}, \quad f \in \Ltwo .
  \end{equation*}

It is natural to choose real spherical harmonics \cite{muller1966spherical}
  \begin{equation*}
    \left\{ \Ylk \,:\, 1 \leq k \leq 2\ell+1,\, 0 \leq \ell \leq t \right\},
  \end{equation*}
as an orthonormal basic for $\Pt$. Noting that the normalisation is such that $ Y_{0,1} = 1/\sqrt{4\pi} $. 
Then
  \begin{equation*}
    \Pt = \mathrm{span} \{ \Ylk : \ell=0,\,\ldots,\,t,\,k=1,\,\ldots,\,2\ell+1 \},
  \end{equation*}
and then the dimension of $\Pt$ is $ d_{t} := \mathrm{dim}(\Pt) = \sum_{\ell=0}^{t} (2\ell+1) = (t+1)^{2} $.
For $t \geq 1$, let the spherical harmonic matrix $\Yt$ be denoted by
  \begin{equation*}
    \Yt(\XN) := [\Ylk(\bm{x}_{j})], \quad k = \sequ{0}{2\ell+1},\;\ell=\sequ{1}{t}; \quad j=\sequ{1}{N}.
  \end{equation*}
It is very important to note the addition theorem of spherical harmonics \cite{atkinson2012spherical}
  \begin{equation} \label{eqt:addthm}
    \sum_{k=1}^{2\ell+1} \Ylk(\mathbf{x})\Ylk(\mathbf{y}) = \frac{2\ell+1}{4\pi} P_{\ell}(\mathbf{x} \cdot \mathbf{y}) \quad \forall \; \mathbf{x},\,\mathbf{y} \in \Stwo,
  \end{equation}
where $ \mathbf{x} \cdot \mathbf{y} $ denotes the usual Euclidean inner product of $\mathbf{x}$ and $\mathbf{y}$ in $\Real^{3}$, and $P_{\ell}$ is the normalized Legendre polynomial of degree $\ell$.
For applications of the addition theorem \eqtref{eqt:addthm}, we refer to \cite{atkinson2012spherical}.

\subsection{Spherical \texorpdfstring{$t$-}{t-}designs }

In \cite{delsarte1977spherical}, lower bounds of even and odd $t$, for the number of nodes $N$ to consist a \std{} are established as following:
  \begin{equation} \label{eq:stdlowbd}
    N(t) \geq
      \begin{cases}
        \frac{ (t+1)(t+3) }{ 4 }, \; & t\,\text{  is odd}, \\
        \frac{  (t+2)^{2}  }{ 4 }, \; & t\,\text{  is even}. \\
      \end{cases}
  \end{equation}
\Std{} achieved this bound \eqtref{eq:stdlowbd} are called tight.
However, Bannai and Damerell \cite{bannai1979tight, bannai1980tight} proved tight \std{} only exists for $t=1,2,3,5$ on $\Stwo$. 
For practical computation, we have to construct \std{} with large $t$.
There is a survey paper on \stds{} given by Bannai and Bannai \cite{bannai2009survey}.
Interval methods are applied to construct reliable computational \stds{} with rigorous proof \cite{an2010well, chen2011computational, chen2006existence}.
In this paper, we are considering two kinds of \std{} as follows:

\subsubsection{Well condition spherical \texorpdfstring{$t$}{t}-designs }
\cite{an2010well} extends the work of \cite{sloan2004extremal} for the case $ N = (t+1)^2 $ by including a constraint that the set of points $\XN$ is a \std{}, as suggested in \cite{chen2006existence}, to extremal \std{} which is a \std{} for which the determinant of a Gram matrix, or equivalently the product of the singular values of a basis matrix, is maximized.
This can be written as the following optimization problem:
\begin{equation} \label{eq:wstdopt}
  \begin{aligned}
    \max_{\XN\subset \Stwo} \;&\; \log\det( \bm{G}_{t}(\XN) ) \\
    \mathrm{s.t.}\;&\;\Ct(\XN) = \bm{E} \bm{G}_{t}(\XN)\bm{e} = \mathbf{0},
  \end{aligned}
\end{equation}
where
\begin{equation*}
  \begin{aligned}
    & \bm{e} := [1,\,\ldots,\,1]^{T} \in \Real^{N}, \\ 
    & \bm{E} := [\bm{e},-\bm{I}_{N-1}] \in \Real^{(N-1)\times N}, \\ 
    & \bm{G}_{t}(\XN) := \Yt(\XN)^{T} \Yt(\XN) \in \Real^{N\times N} .
  \end{aligned}
\end{equation*}
After solving \eqtref{eq:wstdopt} by nonlinear optimization methods, the interval analysis provides a series of narrow intervals, which contain computational \std{} and a true \std{}.
Consequently, the mid point of these intervals are determinated as well conditioned \std{}, for detail, see \cite{an2010well}. \par

Following the methods in \cite{an2010well}, we use extremal systems \cite{sloan2004extremal} which maximize the determinant without any additional constraints as the starting points to solve this problems \eqtref{eq:wstdopt}.
We also use interval methods, which memory usage is optimized, to prove that close to the computed extremal \std{} there are exact \std{}.
Finally, we obtain well conditioned \std{} with degree up to 160.
For detail, we refer to another paper on construct well conditioned \std{} for $t$ up to 160, see \cite{an2016numerical}.

\subsubsection{Efficient spherical \texorpdfstring{$t$}{t}-designs }

In \cite{womersley2015efficient}, Womersley introduced efficient \stds{} with roughly $\frac{1}{2}t^{2}$ points.
The point number is close to the number in a conjecture by Hardin and Sloane that $ N=\frac{1}{2}t^{2}+o(t^2) $ \cite{hardin1996mclaren}.
The author used Levenberg-Marquardt method to solve the following problem
\begin{equation}
  \begin{aligned}
    \min_{\XN\subset \Stwo} \;&\; A_{N,t}(\XN) = \dfrac{4\pi}{N^{2}} \sum_{\ell=1}^{t} \sum_{k=1}^{2\ell+1} \left( \sum_{j=1}^{N} \Ylk(\mathbf{x}_{j}) \right)^{2} .
  \end{aligned}
\end{equation}
This point sets can be download at \url{http://web.maths.unsw.edu.au/~rsw/Sphere/EffSphDes/index.html}.

%% file: chapter/3-wce.tex
\section{Worst-case error of spherical \texorpdfstring{$t$}{t}-designs }
\label{sec:wce}

This section considers the worst-case error for numerical integration over $\Stwo$ \cite{brauchart2014qmc} \cite{hesse2005worst}. 
In this section we follow notations and definitions from \cite{brauchart2014qmc}.
The Sobolev space, denoted by $ \Hs := \Hs(\Stwo)$, can be defined for $s \geq 0$ as the set of all functions $ f \in \Ltwo $ with Laplace-Fourier coefficients
  \begin{equation*}
    \hat{f}_{\ell,k} := \inprod{f}{\Ylk}_{\Ltwo} = \int_{\Stwo} f(\mathbf{x}) \Ylk(\mathbf{x}) \domegax,
  \end{equation*}
satisfying
  \begin{equation*}
    \sum_{\ell=0}^{\infty} \sum_{k=1}^{2\ell+1} (1+\lambda_{\ell})^{s} \abs{ \hat{f}_{\ell,k} }^2 < \infty,
  \end{equation*}
where $ \lambda_{\ell} := \ell(\ell+1) $. Obviously, by letting $ s=0 $, then $ \mathbb{H}^{0}= \Ltwo $. 
The norm of $ \Hs $ can be defined as
  \begin{equation*}
    \norm{f}{\Hs} = \left[ \sum_{\ell=0}^{\infty} \sum_{k=1}^{2\ell+1}  \frac{1}{a_{\ell}^{(s)}} \hat{f}_{\ell,k}^2 \right]^{\frac{1}{2}},
  \end{equation*}
where the sequence of positive parameters $\alpha_{\ell}^{(s)}$ satisfies $a_{\ell}^{(s)} \sim (1+\lambda_{\ell})^{-s} \sim (\ell+1)^{-2s}$.

The \emph{worst-case error} of the \std{} $\XN$ on $\Hs$ can be defined as
  \begin{equation}\label{eqt:wce_std}
    \wce(Q[\XN]) := \sup\left\{ \abs{Q[\XN](f)-I(f)} \,\big\vert\, f\in\Hs,\, \norm{f}{\Hs} \leq 1  \right\},
  \end{equation}
where $Q[\XN](f) := \frac{4\pi}{N} \sum_{j=1}^{N} f(\mathbf{x}_{j}) $. \par

Before introducing the formula of worst-case error, we show the signed power of the distance, with sign $(-1)^{L+1}$ with $ L := L(s) := \lfloor s-1 \rfloor $,  that has the following Laplace-Fourier expansion \cite{atkinson2012spherical}: for $\mathbf{x},\,\mathbf{y}\in\Stwo$,
  \begin{equation*}
    (-1)^{L+1}\abs{\mathbf{x}-\mathbf{y}}^{2s-2} = (-1)^{L+1}V_{2-2s}(\Stwo) + \sum_{\ell=1}^{\infty} \alpha_{\ell}^{(s)} Z(2, \ell) P_{\ell}(\mathbf{x} \cdot \mathbf{y}),
  \end{equation*}
where $P_{\ell}$ is the normalized Legendre polynomial,
\begin{align*}
    & V_{2-2s}(\Stwo) = \int_{\Stwo}\int_{\Stwo} \abs{\mathbf{x}-\mathbf{y}}^{2s-2} \domegax \df\omega(\mathbf{y}) = 2^{2s-1} \frac{ \Gamma{(3/2)} \Gamma{(s)} }{ \sqrt{\pi}\Gamma{(1+s)} } = \frac{ 2^{2s-2} }{ s }, \\
    & \alpha_{\ell}^{(s)} := V_{2-2s}(\Stwo) \frac{ (-1)^{L+1} (1-s)_{\ell} }{ (1+s)_{\ell} }
      = V_{2-2s}(\Stwo) \frac{ (-1)^{L+1} \Gamma{(1-s+\ell)}\Gamma{(1+s)} }{ \Gamma{(1+s+\ell)}\Gamma{(1-s)} }, \\
    & Z(2, \ell) = (2\ell+2-1) \frac{ \Gamma{(\ell+2-1)} }{ \Gamma{(2)}\Gamma{(\ell+1)} } =  2\ell+1  .
\end{align*}
From \cite{brauchart2014qmc}, we know that worst-case error is divided into two cases:
\begin{enumerate}[label={\bf Case \Roman* }, ref={\bf Case \Roman* }]
  \item For $ 1 < s < 2 $ and $L=0$, the worst-case error is given by
    \begin{equation}
      \wce(Q[\XN]) = \left( V_{2-2s}(\Stwo) - \frac{1}{N^2}\sum_{j=1}^{N}\sum_{i=1}^{N} \abs{\mathbf{x}_{j}-\mathbf{x}_{i}}^{2s-2} \right)^{\frac{1}{2}}.
    \end{equation}
  \item For $ s > 2 $ and $L$ satify $L:=L(s)=\lfloor s-1  \rfloor$, the worst-case error is given by
    \begin{equation}
      \begin{aligned}
        \wce(Q[\XN]) = \Bigg( \frac{1}{N^2}\sum_{j=1}^{N}\sum_{i=1}^{N} \Big[ \mathcal{Q}_{L}(\mathbf{x}_{j}\cdot \mathbf{x}_{i}) + (-1)^{L+1}\abs{\mathbf{x}_{j}-\mathbf{x}_{i}}^{2s-2} \Big] \\
          -(-1)^{L+1}V_{2-2s}(\Stwo) \Bigg)^{\frac{1}{2}},
      \end{aligned}
    \end{equation}
    where 
      \begin{equation*}
        \mathcal{Q}_{L}(\mathbf{x}_{j}\cdot \mathbf{x}_{i}) = \sum_{\ell=1}^{L} \left( (-1)^{L+1-\ell}-1 \right) \alpha_{\ell}^{(s)} Z(2, \ell) P_{\ell}(\mathbf{x}_{j}\cdot \mathbf{x}_{i}). 
      \end{equation*}
\end{enumerate}

\subsection{Numerical experiments on worst-case error}
By using the definition of worst-case error, we calculate and compare worst-case error of two \stds{}: well condition \std{}\cite{an2010well} and efficient \std{}\cite{womersley2015efficient}.
\figref{fig:wce} gives, when $s=1.5,\, 2.5,\, 3.5\; \text{and}\; 4.5 $, worst-case error for WSTD and ESTD \cite{womersley2015efficient}.
Worst-case error for WSTD is smaller which means that it has a better performance in numerical integration when the precision $t$ of two point sets are the same.

\begin{figure}
  \centering
  \subfigure[$s=1.5$ \label{fig:wce_3_2}]{
    \includegraphics[width=0.7\textwidth]{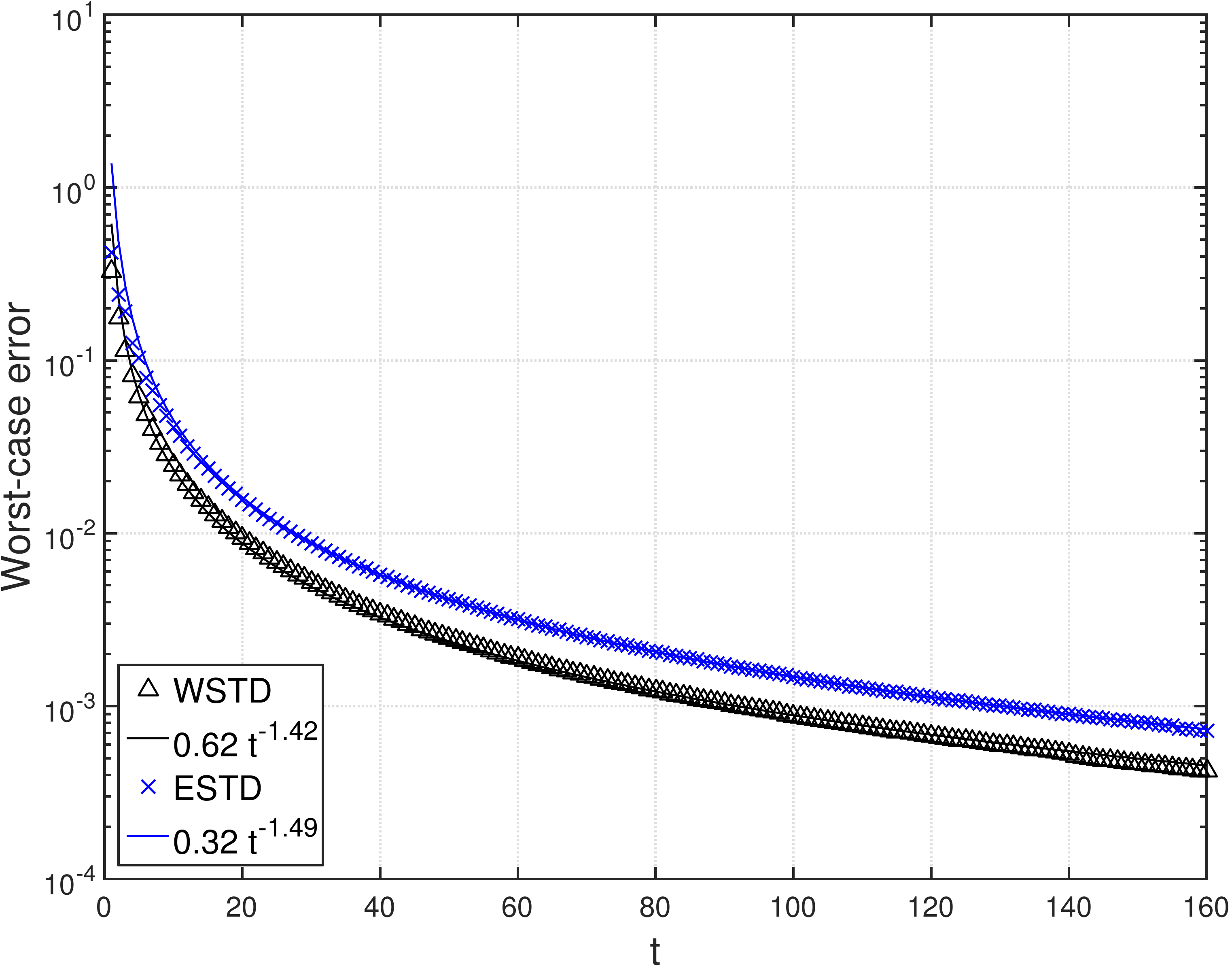}
  }
  \subfigure[$s=2.5$ \label{fig:wce_5_2}]{
    \includegraphics[width=0.7\textwidth]{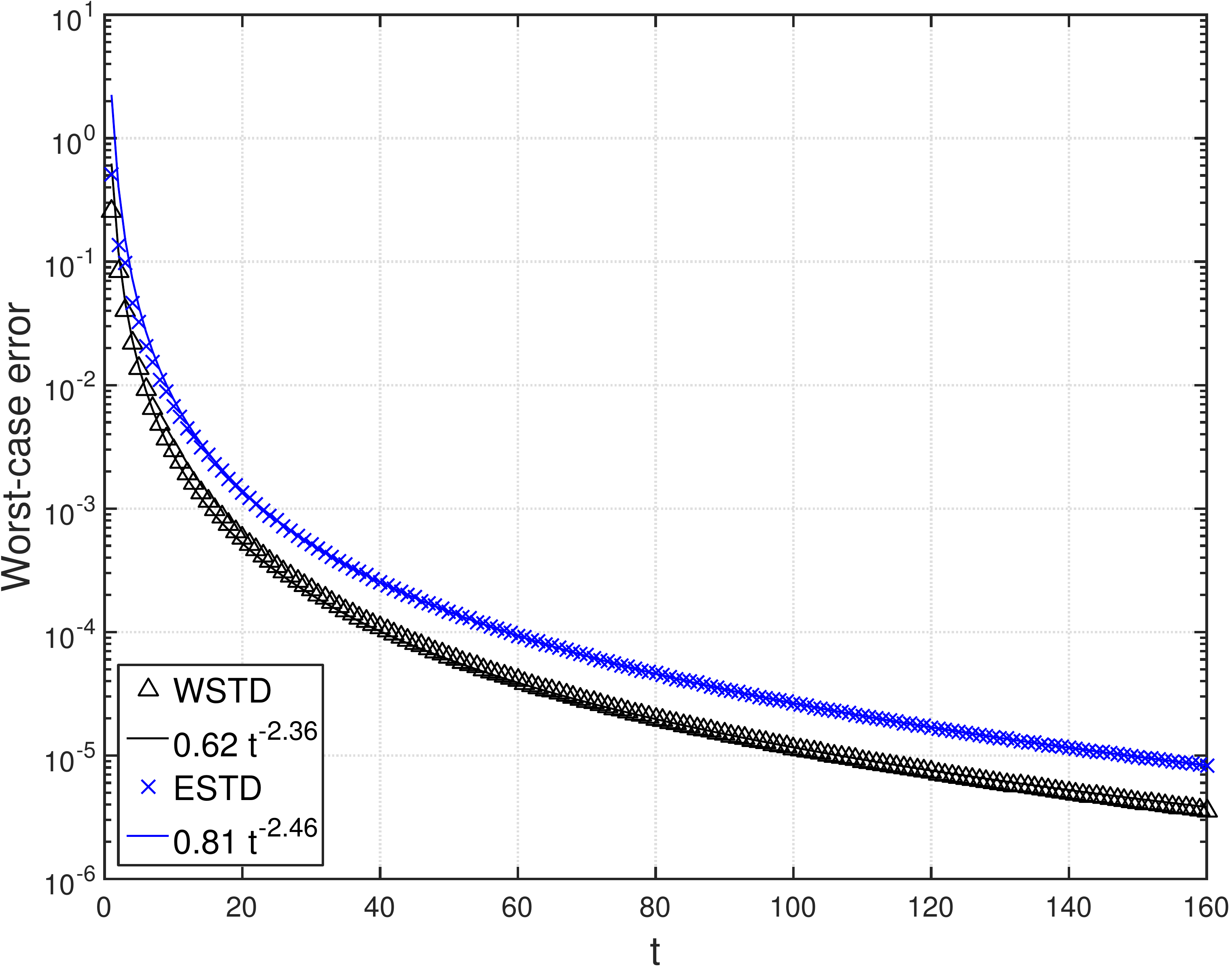}
  }
  \caption{Worst-case error of two point sets }\label{fig:wce}
\end{figure}
\begin{figure}
  \centering
  \contsubfigure[$s=3.5$  \label{fig:wce_7_2} ]{
    \includegraphics[width=0.7\textwidth]{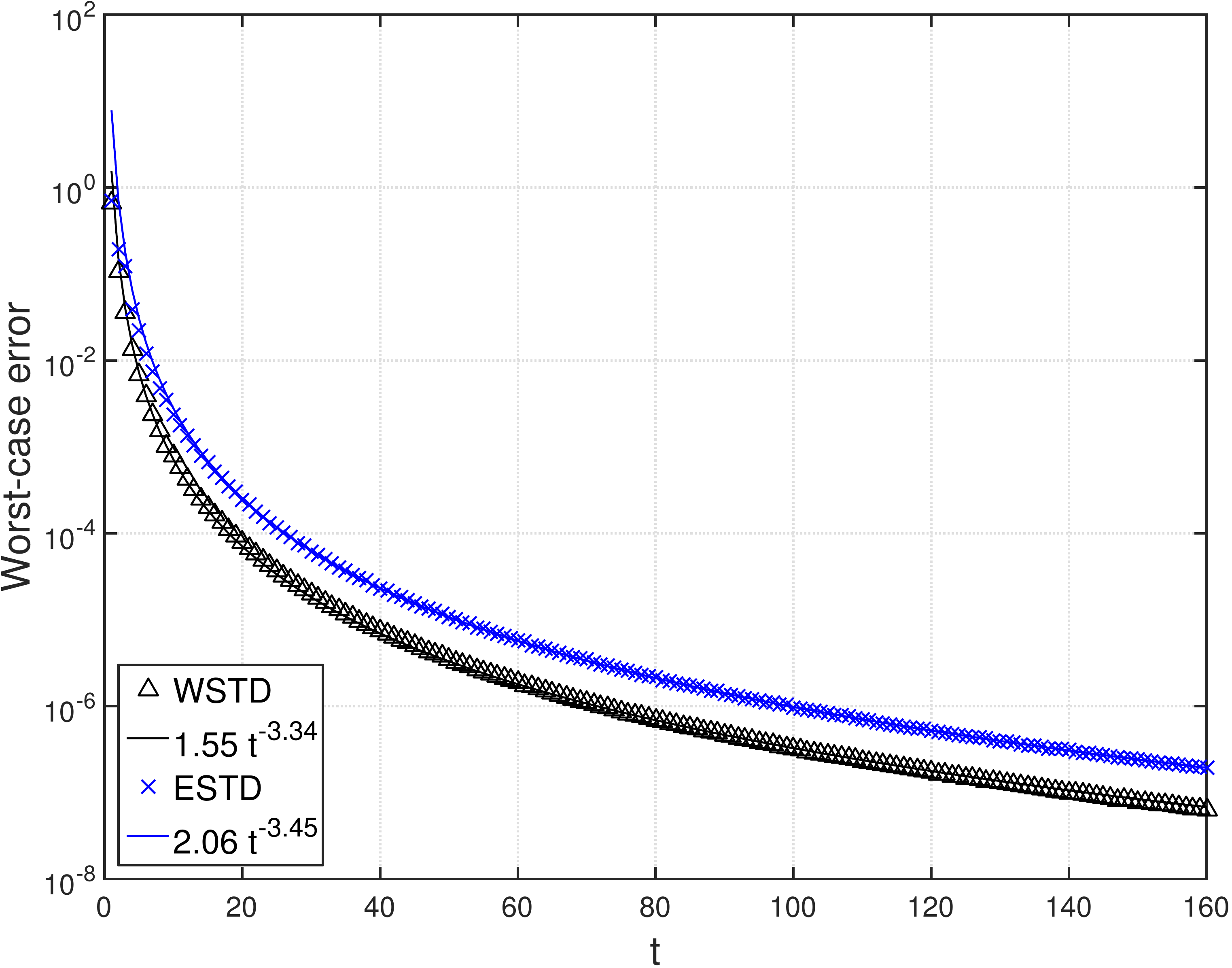}
  }
  \contsubfigure[$s=4.5$ \label{fig:wce_9_2}]{
    \includegraphics[width=0.7\textwidth]{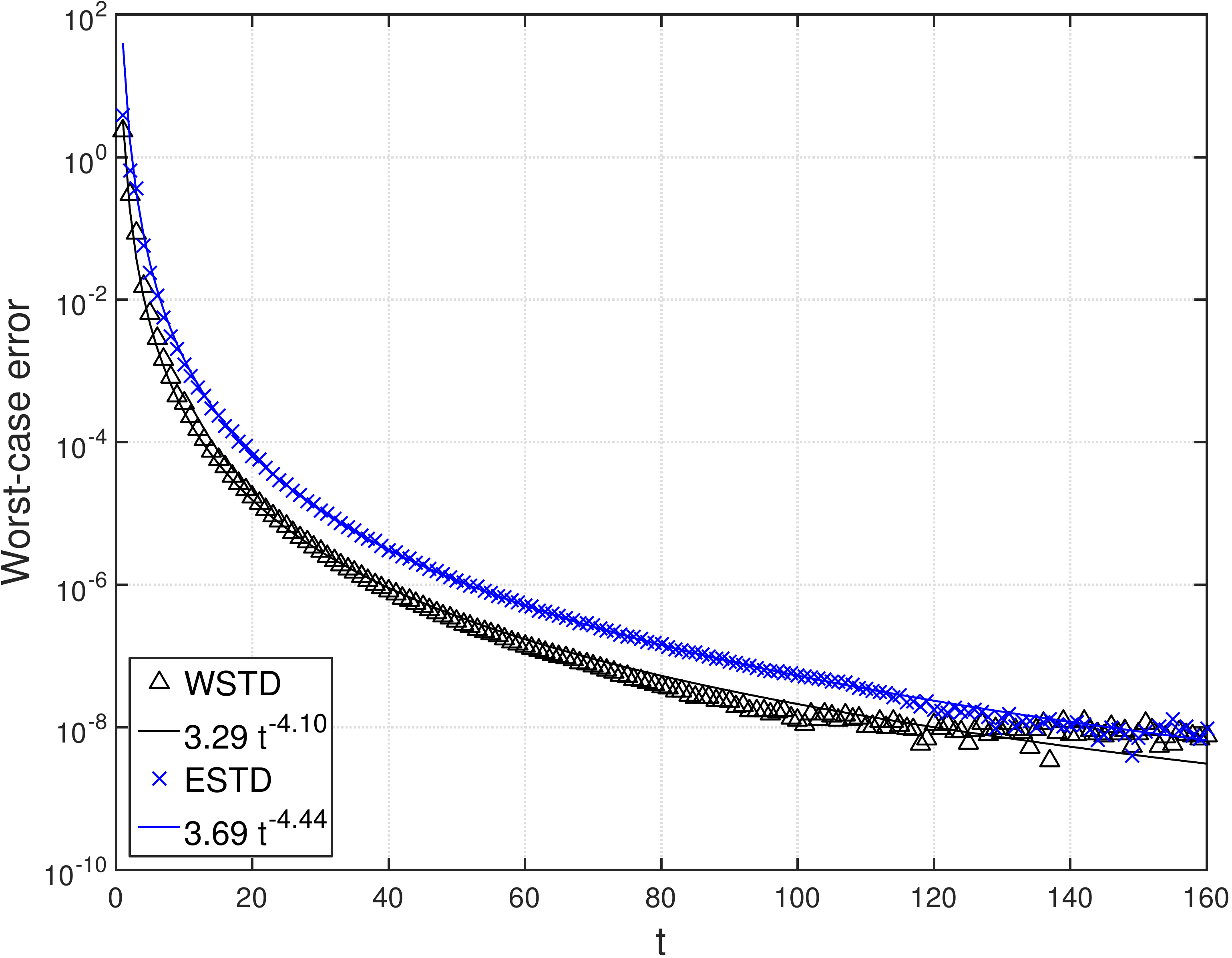}
  }
  \contcaption{Worst-case error of two point sets }
  \subconcluded
\end{figure}

\subsection{Conjecture on worst-case error of spherical \texorpdfstring{$t$}{t}-design}
From the above interesting numerical experiments, we propose a reasonable conjecture as following:
\begin{conjecture}
  Let  $\wce(Q[\XN])$ be worst-case error of the \std{}(see \eqtref{eqt:wce_std}). Then when $N$ increases, $\wce(Q[\XN])$ decreases. That is to say:
  \begin{equation*}
    \wce(Q[\mathcal{X}_{N'}])\leq \wce(Q[\XN]),\quad N'>N \ \text{for fixed}\ t.
  \end{equation*}
\end{conjecture}

%% file: chapter/4-tran.tex
\section{Transformations for singular functions on \texorpdfstring{$\Stwo$}{S2} }
\label{sec:tran}

In this section, we consider two variable transformations for the approximation of spherical integral $\label{eqt:singular_integral} I(f) $ in which $ f(\mathbf{x}) $ is singular at a point $ \mathbf{x}_{0} $.
Examples are the single layer integral
  \begin{equation*}
    \dpsty \int_{\Stwo} \dfrac{ g(\mathbf{x}) }{ \abs{\mathbf{x}-\mathbf{x}_{0}} }  \domegax
  \end{equation*}
and the double layer integral
  \begin{equation*}
    \int_{\Stwo} \dfrac{ g(\mathbf{x}) \abs{ (\mathbf{x}-\mathbf{x}_{0})\cdot \bm{n}_{\mathbf{x}} } }{ \abs{\mathbf{x}-\mathbf{x}_{0}}^{3} }  \domegax ,
  \end{equation*}
where $g(\mathbf{x})$ is smooth function, $ \bm{n}_{\bm{x}} $ is the outward normal to $\Stwo$ at $\mathbf{x}$, and $ (\mathbf{x}-\mathbf{x}_{0})\cdot \bm{n}_{\mathbf{x}} $ is the dot product of two vectors $(\mathbf{x}-\mathbf{x}_{0})$ and $\bm{n}_{\mathbf{x}}$.
From \cite{sidi2005application}, we know the double-layer integral is simply $1/2$ times the single-layer integral over $\Stwo$.
So it is enough to treat the single-layer case.
In this case of that integrand is a singular function, we use a variable transformation, such as Atkinson's transformation \cite{atkinson2004quadrature} and Sidi's transformation \cite{sidi2005application}, rather than approximating this integral directly.
In the following, we introduce two variable transformations : Atkinson's transformation and Sidi's transformation.

\subsection{Atkinson's transformation \texorpdfstring{$\Tone$}{T1} }

We consider the transformation $ \Tone: \Stwo \xlongrightarrow[\text{onto}]{1-1} \Stwo $ introduced by Atkinson \cite{atkinson2004quadrature}.
We show $\tilde{\mathbf{x}} \in \Stwo$  by standard spherical coordinates with $\theta$ (polar angle, $0\leq\theta<\pi$) and $\phi$ (azimuth angle, $0\leq\theta<2\pi$). Define the transformation:
  \begin{equation*}
    \Tone: \mathbf{x}=(\cos\phi\sin\theta ,\, \sin\phi\sin\theta ,\, \cos\theta)^{T}
      \mapsto \tilde{\mathbf{x}} = \frac{ \left( \cos\phi\sin^{q}\theta ,\, \sin\phi\sin^{q}\theta ,\, \cos\theta \right)^{T} }{ \sqrt{\cos^{2}\theta+\sin^{2q}\theta} } . 
  \end{equation*}
In this transformation, $ q\geq 1$ is a `grading parameter'. Then north and south poles of $\Stwo$ remain fixed, while the region around them is distorted by the mapping. Then the integral $I(f)$ becomes
  \begin{equation*}
    I(f) = \int_{\Stwo} f( \Tone(\mathbf{x}) ) J_{\Tone}(\mathbf{x}) \domegax,
  \end{equation*}
with $J_{\Tone}(\mathbf{x})$  the Jacobian of the mapping $\Tone$,
  \begin{equation*}
    J_{\Tone}(\mathbf{x}) = \frac{ \sin^{2q-1}\theta (q\cos^{2}\theta+\sin^{2}\theta) }{ (\sin^{2q}\theta+\cos^{2}\theta)^{\frac{3}{2}} }, \quad q \geq 1.
  \end{equation*}
As shown in \cite{atkinson2004quadrature, atkinson2005quadrature, sidi2005analysis}, for smooth integrand, when $2q$ is an odd integer, trapezoidal rule enjoys a fast convergence.
Consequently, in the following numerical experiments, we set grading parameter $q$ such  that $2q$ is an odd integer.

\subsection{Sidi's transformation \texorpdfstring{$\Ttwo$}{T2} }
Sidi \cite{sidi2005application} introduced another variable transformation $\Ttwo$ with the aid of spherical coordinate as following:
  \begin{equation*}
    \Ttwo: \mathbf{x}=(\theta, \phi)^{T}
    \mapsto \tilde{\mathbf{x}} = \left(\Psi \left( {\textstyle\frac{\theta}{2\pi}} \right), \phi \right)^{T},
  \end{equation*}
where $\Psi(t)$ is derived from a standard variable transformation $\psi(t)$ in the extended class $\Ttwo$ of Sidi \cite{sidi2006extension}, and $\Psi(t) = \pi \psi(t)$, which is the first way to do variable transformation in \cite{sidi2005application}.
The standard variable transformation $\psi(t)$, just as the original $\sin^{m}$ -- transformation, is defined via
  \begin{equation}\label{eqt:psi_m}
    \psi_{m}(t) = \frac{\Theta_{m}(t)}{\Theta_{m}(1)}; \quad \Theta_{m}(t) = \int_{0}^{t} (\sin\,\pi u)^{m} \df u.
  \end{equation}
Here $m \geq 1$ act as the `grading parameter' in $\Tone$.
From $\Theta_{m}(t)$'s derivative $ \Theta_{m}'(t) = (\sin{\pi t})^m $, we have $ \psi_{m}'(t)= (\sin{\pi t})^m / \Theta_{m}(1) $.
Obviously, $\Theta_{m}'(t)$ is symmetric with respect to $ t=1/2 $. 
So $\Theta_{m}(t)$ satisfies equation $ \Theta_{m}(t) = \Theta_{m}(1) - \Theta_{m}(1-t) $ for $ t \in \left[ 1/2,\, 1 \right] $. 
Thus, $ \Theta_{m}(1) = 2\,\Theta_{m}(1/2) $. 
Consequently,
  \begin{equation}\label{eqt:psi_symm}
    \psi_{m}(t) = \frac{\Theta_{m}(t)}{2\Theta_{m}(1/2)} \;\; \text{for} \;\; t \in  \left[ 0,\, 1/2 \right]; \quad \psi_{m}(t) = 1-\psi_{m}(1-t)  \;\; \text{for}  \;\;  t \in  \left[ 1/2,\, 1 \right].
  \end{equation}
From equality
  \begin{equation*}
    \Theta_{m}(t) = \frac{m-1}{m} \Theta_{m-2}(t) - \frac{1}{\pi m} (\sin{\, \pi t})^{m-1} \cos{\, \pi t},
  \end{equation*}
we have the recursion relation
  \begin{equation} \label{eq:psirecur}
    \psi_{m}(t) = \psi_{m-2}(t) - \frac{ \Gamma{(m/2)} }{ 2\sqrt{\pi} \Gamma{((m+1)/2)} } \, (\sin{\, \pi t})^{m-1} \cos{\, \pi t}.
  \end{equation}
When $m$ is a positive integer, $ \psi_{m}(t) $ can be expressed in terms of elementary functions.
In this case, $ \psi_{m}(t) $ can be computed via the recursion relation \eqref{eq:psirecur}, with the initial conditions
  \begin{equation}
    \psi_{m}(t) = t \quad \text{and} \quad \psi_{m}(t) = \frac{1}{2}\left( 1-\cos{\pi t} \right).
  \end{equation}
When $m$ is not an integer,\ $ \psi_{m}(t) $ cannot be expressed in terms of elementary functions.
However, it can be expressed conveniently in terms of hypergeometric functions.
Because of symmetry and \eqtref{eqt:psi_symm}, it is enough to consider the computation of $ \Theta_{m}(t) $ only for $ t \in  \left[ 0,\, 1/2 \right] $. One of the representations in terms of hypergeometric functions now reads
  \begin{equation}\label{eqt:theta_2}
    \Theta_{m}(t) = \frac{(2K)^{m+1}}{\pi(m+1)} F \left( \frac{1}{2}-\frac{1}{2}m,\, \frac{1}{2}+\frac{1}{2}m;\, \frac{1}{2}m+\frac{3}{2};\, K^{2} \right); \quad K=\sin{\frac{\pi t}{2}}.
  \end{equation}
Then the expression of $ \psi_{m}(t) $ in \eqtref{eqt:psi_m} follows from \eqtref{eqt:theta_2}. \par
Now the integral $I(f)$ becomes
  \begin{equation*}
    I(f) = \int_{\Stwo} f( \Ttwo(\mathbf{x}) ) J_{\Ttwo}(\mathbf{x}) \domegax,
  \end{equation*}
with Jacobian 
  \begin{equation*}
    J_{\Ttwo}(\mathbf{x}) = \Psi_{m}'\left( {\textstyle\frac{\theta}{2\pi}} \right) = {\textstyle\frac{1}{2}} \psi_{m}'\left( {\textstyle\frac{\theta}{2\pi}} \right).
  \end{equation*}

\subsection{Numerical integration method with orthogonal transformation}

Since the above variable transformations are based on that the singular point is at the north pole of $\Stwo$, we need to move the north pole to the singular point.
The original coordinate system of $\Real^{3}$ needs to be rotated to have the north pole of $\Stwo$ in the rotated system be the location of the singularity in integrand.
Atkinson \cite{atkinson2004quadrature} used an orthogonal Householder transformation.
Here we introduce the rotation transformation\ $\bm{R}$\ in $\Real^{3}$,
  \begin{equation}\label{matr:Rot}
    \hat{\mathbf{x}} = \bm{R} \mathbf{x},  \;\; \quad \hat{\mathbf{x}},\,\mathbf{x} \in \Stwo .
  \end{equation}
In fact, $\bm{R}$ can be expressed as follows:
  \begin{equation*}
    \bm{R} = \bm{R}_{z}(\phi)\bm{R}_{y}(\theta) =
    \left[ \begin{array}{ccc} \cos{\phi} & -\sin{\phi} & 0 \\ \sin{\phi} & \cos{\phi} & 0 \\ 0 & 0 & 1 \end{array} \right]
    \cdot
    \left[ \begin{array}{ccc} \cos{\theta} & 0 & \sin{\theta} \\ 0 & 1 & 0  \\ -\sin{\theta} & 0 & \cos{\theta} \end{array} \right],
  \end{equation*}
such that
  \begin{equation*}
    \bm{R}\left[ \begin{array}{c} 0 \\ 0 \\ 1 \end{array} \right] =  \mathbf{x}_{0}.
  \end{equation*}
Here, $\mathbf{x}_{0} \in \Stwo$ is the singular point. Obviously, when the singular point is just at the north pole, the rotation matrix\ $\bm{R}$\ is the identity matrix\ $\bm{I}$.
\par

By using \eqtref{eqt:quadrarule}, the singular  integral can by approximated by the following form with these positive weight quadrature rules:
  \begin{equation*}
    I(f) \approx \sum_{j=1}^{N} \omega_{j} f( \bm{R}\mathcal{T}_{i}(\mathbf{x}_{j}) ) J_{\mathcal{T}_{i}}(\mathbf{x}_{j}), \quad i=1,2.
  \end{equation*}
Here, $\Tone$, $\Ttwo$ correspond to Atkinson's transformation and Sidi's transformation, respectively.

%% file: chapter/5-numres.tex
\section{Numerical Results}\label{sec:num_res}

\subsection{Quadrature nodes over \texorpdfstring{$\mathbb{S}^2$}{S2} }
\begin{figure}
  \centering
  \subfigure[Bi. trapezoidal rule, $N=231$ \label{fig:tr}]{
    \includegraphics[width=0.30\textwidth]{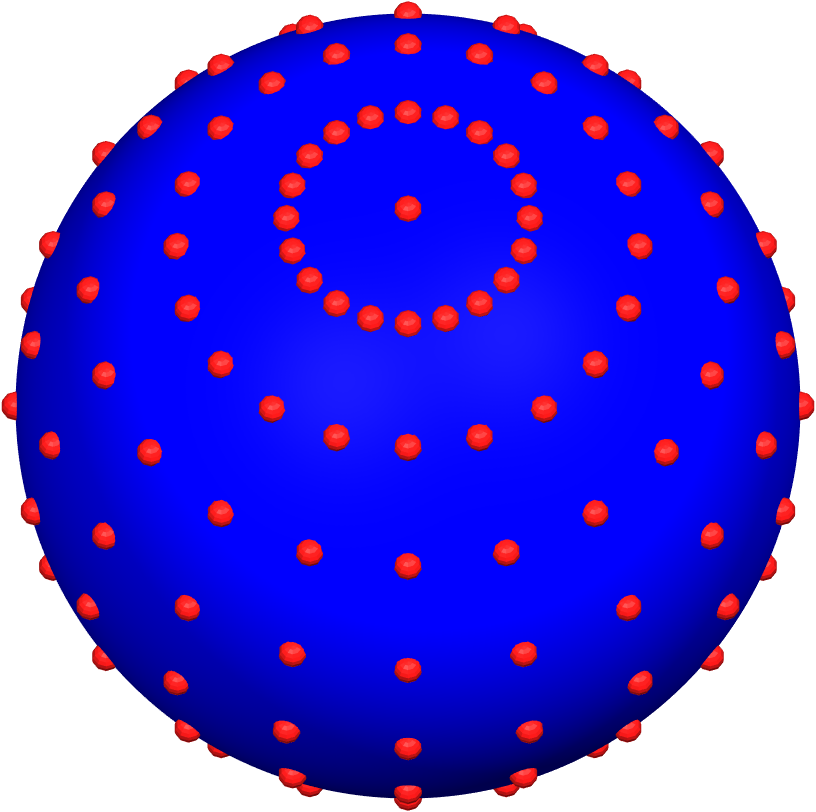}
  }
  \subfigure[Equal area points, $N=225$ \label{fig:eq}]{
    \includegraphics[width=0.30\textwidth]{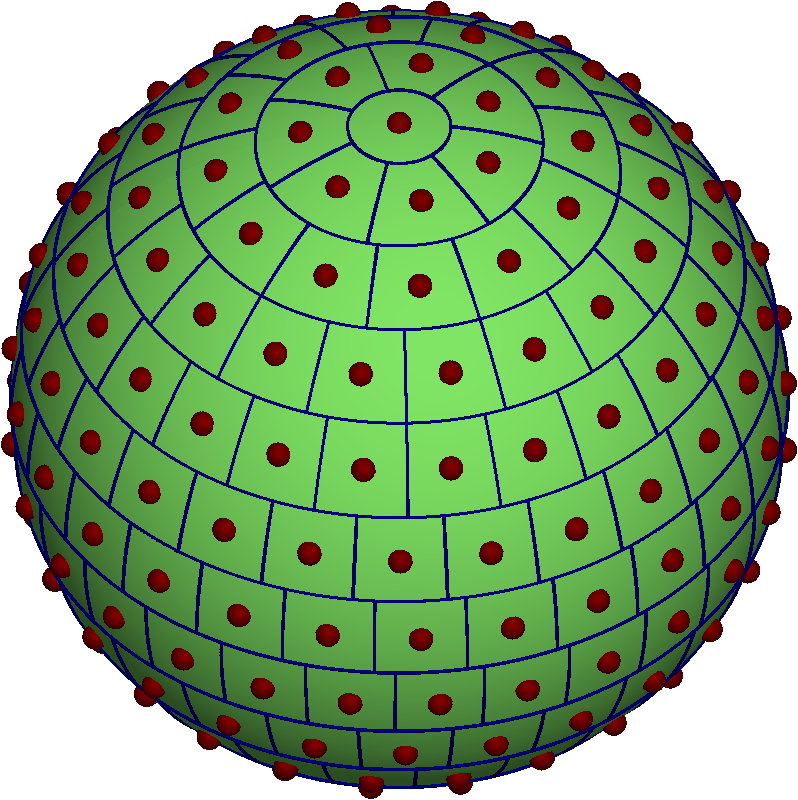}
  }
  \subfigure[\std{}, $N=225$ \label{fig:std}]{
    \includegraphics[width=0.30\textwidth]{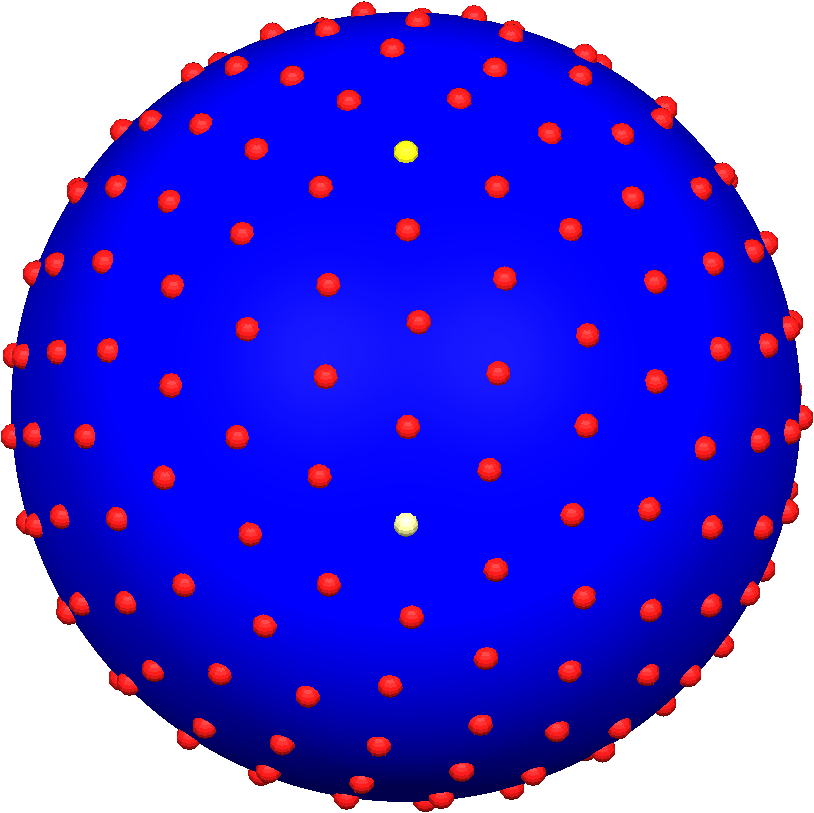}
  }
  \caption{Point sets on sphere}\label{fig:pointset}
  \subconcluded
\end{figure}

In this section we investigate three quadrature rules to approximate the integration of serval test functions:
\begin{itemize}
  \item {\bf Bivariate trapezoidal rule } This quadrature rule is consist of Longitude-Latitude points which divide the longitude and latitude equally. By using the spherical coordinate $(\theta,\,\phi)$, for $n \geq 1$, let $ h=\pi/n $, and $ \theta_{j} = \phi_{j} = jh $.
      Bivariate trapezoidal can be written in the following formula
      \begin{equation*}
        \QN(f) := \frac{\pi^2}{n^2} \sum_{j=0}^{2n}{}^{''} \sum_{i=0}^{n}{}^{''} f(\theta_{i}, \phi_{j}).
      \end{equation*}
      Here  the superscript notation ${}{''}$ means to multiply the first and last terms by $ 1/2 $ before summing. Atkinson \cite{atkinson2004quadrature} and Sidi \cite{sidi2005application} added a transformation to led rapid convergence or reduce the effect of any singularities in $f$ at the poles.
      In the following numerical experiment, for continuous function, we use Atkinson's transformation with $q=2.5$, which shows faster convergence than $2$ and $3$, and for singularity we use another $q$.
  \item {\bf Equal Partition Area points on sphere } This integration rule based on partitioning the sphere into a set of $N$ open domains $T_{j} \subset \Stwo,\,j=\sequ{1}{N}$, that is $ T_{j} \bigcap T_{k} = \emptyset $ for $ j \neq k$, and $\dpsty \bigcup_{j=1}^{N} \overline{T_{j}}=\Stwo$ where $\overline{T_{j}}$ represent the closure of $T$.
      So the quadrature weight is the surface area of  $T$.
      In \cite{leopardi2009diameter},  Leopardi developed an algorithm to divide the sphere into $N$ equal area partition efficiently.
      The center of each partition is chosen as the quadrature node.
      So the corresponding quadrature rule is of equal positive weight:
        \begin{equation*}
          \QN(f) := \frac{4\pi}{N} \sum_{i=1}^{N} f(\mathbf{x}_{i}).
        \end{equation*}
  \item {\bf Well conditional \std{}}\ ( $t$ up to 160 ).
\end{itemize}
The direct observation of these three point sets can be found in \figref{fig:pointset}.

\subsection{Geometric properties}
In this section we concentrate on the geometric properties of point sets over sphere.
Naturally, the distance between any two points $\mathbf{x}$ and $\mathbf{y}$ on $\Stwo$ is measured by the geodesic distance:
  \begin{equation*}
    \dist(\mathbf{x},\, \mathbf{y})=\cos^{-1}(\mathbf{x}\cdot \mathbf{y}), \quad \mathbf{x},\,\mathbf{y}\in\Stwo,
  \end{equation*}
which is the natural metric on $\Stwo$.
The quality of the geometric distribution of point set $\XN$ is often characterized by the following two quantities and their ration. 
The \emph{mesh norm}
  \begin{equation} \label{eqt:meshnorm}
    h_{\XN} := \max_{\mathbf{y}\in\Stwo} \min_{\mathbf{x}_{i}\in\XN} \dist(\mathbf{y}, \mathbf{x}_{i})
  \end{equation}
and the \emph{minimal angle}
  \begin{equation*}
    \delta_{\XN}=\min_{ \mathbf{x}_{i},\mathbf{x}_{j}\in\XN,i\neq j } \dist(\mathbf{x}_{i}, \mathbf{x}_{j}).
  \end{equation*}
The mesh norm is the covering radius for covering the sphere with spherical caps of the smallest possible equal radius centered at the points in $\XN$ , while the minimal angle $\delta_{\XN}$ is twice the packing radius, so $h_{\XN}\geq\delta_{\XN}/2$.
The \emph{mesh ratio} $\rho_{\XN}$
  \begin{equation*}
    \rho_{\XN}:=\frac{ 2h_{\XN} }{ \delta_{\XN} }\geq 1
  \end{equation*}
is a good measure for the quality of the geometric distribution of $\XN$ : the smaller $\rho_{\XN}$ is; the more uniformly are the points distributed on $\Stwo$ \cite{hesse2010numerical}.
\par

The geometric properties of above three point sets are shown in \figref{fig:ptprop}.
Comparing the mesh norm, minimal angle and mesh ratio of three point sets in three subfigures, it can be seen that the mesh norm of WSTD is between the Bivariate trapezoidal rule and Equal area points.
\begin{figure}[htbp]
  \centering
  \subfigure[Mesh norm of three point sets \label{fig:meshnorm3point}]{
    \includegraphics[width=0.65\textwidth]{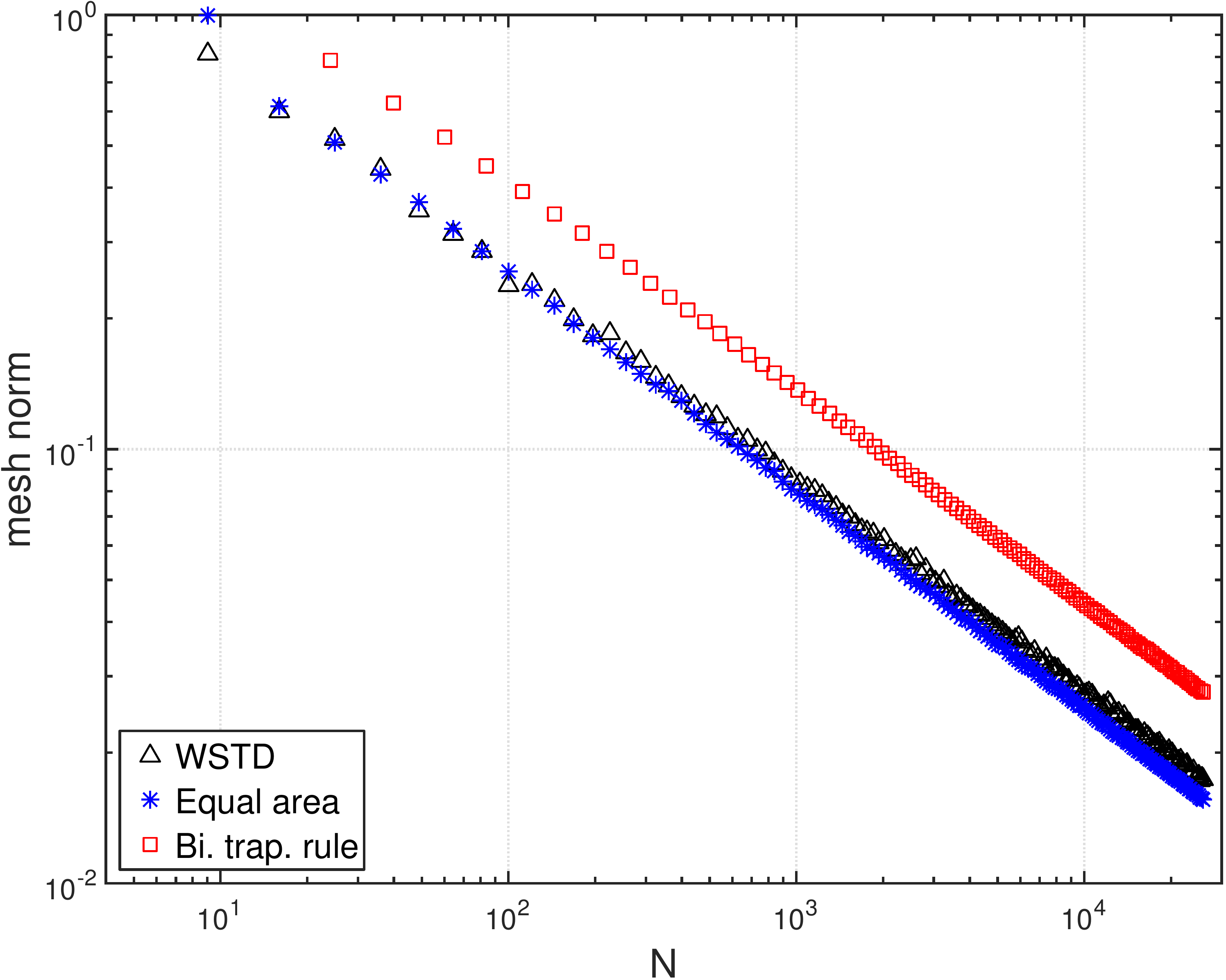}
  }

  \subfigure[Minimal angle of three point sets \label{fig:minangle3point}]{
    \includegraphics[width=0.65\textwidth]{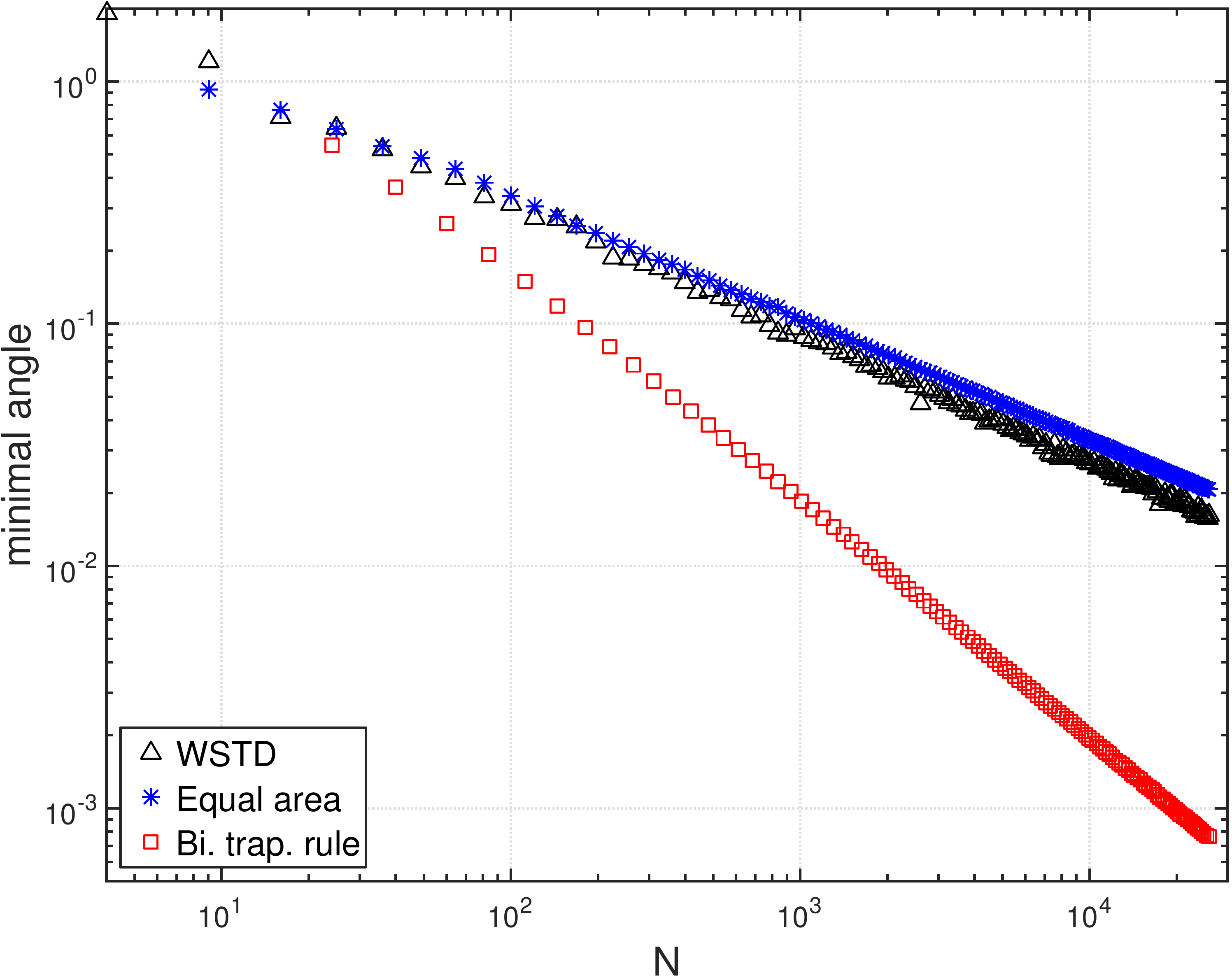}
  }
  \caption{Geometry of above three point sets}\label{fig:ptprop}
\end{figure}
\begin{figure}[htb]
  \centering
  \contsubfigure[Mesh ratio of three point sets \label{fig:meshratio3point}]{
    \includegraphics[width=0.65\textwidth]{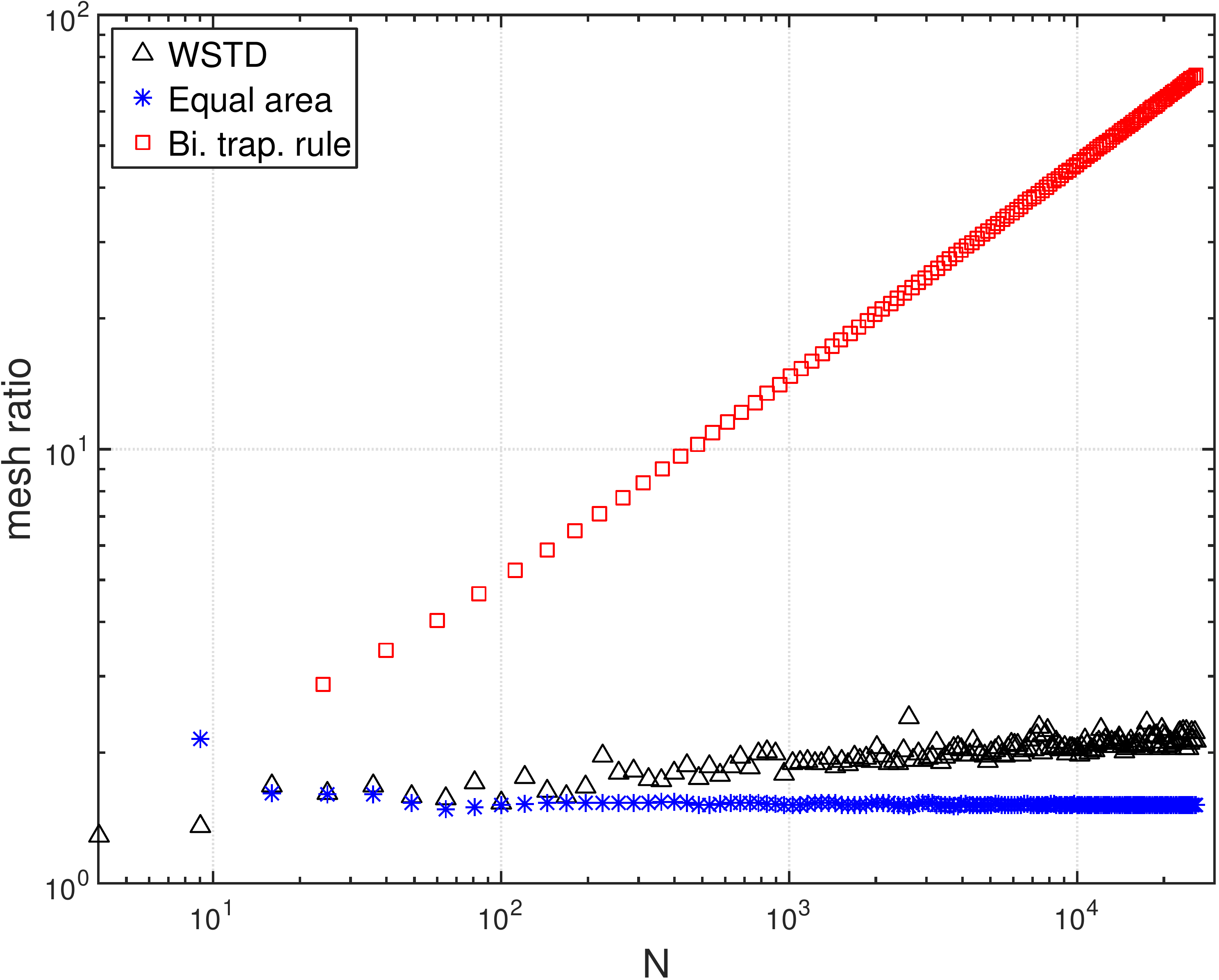}
  }
  \contcaption{Geometry of three point sets}
  \subconcluded
\end{figure}

\subsection{Test functions}

The used functions are expressed as follows.
\begin{align}
&
    \begin{aligned}f_1(x, y, z) = & 0.75 \exp(-(9x-2)^2/4  - (9y-2)^2/4 - (9z-2)^2/4) \\
                                  &+0.75 \exp(-(9x+1)^2/49 - (9y+1)/10  - (9z+1)/10 ) \\
                                  &+0.5  \exp(-(9x-7)^2/4  - (9y-3)^2/4  - (9z-5)^2/4 ) \\
                                  &-0.2  \exp(-(9x-4)^2    - (9y-7)^2   - (9z-5)^2 ),
    \end{aligned}\\
& f_2(x,y,z) = \frac{\sin^2(1+\abs{x+y+z})}{10}, \\
& f_3(x,y,z) = \frac{1}{101-100z}, \\
& f_4(x,y,z) =
  \begin{cases}
   h_0 \cos^2{ \left( \frac{\pi}{2} \cdot \frac{r}{R} \right) } & \textrm{if} \;\; r < R, \\
   0                                                              & \textrm{if}\;\; r \geq R,
  \end{cases}
  \quad r = \dist \big( (x,\,y,\,z)^{T} ,\, (x_c,\,y_c,\,z_c)^{T} \big), 
  \label{eq:capfun}\\
& f_5(x,y,z) = \frac{\exp{(x+2y+3z)}}{\vectornorm{(x,\,y,\,z)^{T}-(x_0,\,y_0,\,z_0)^{T}}}, \\
& f_6(x,y,z) = \frac{\exp{0.1(x+2y+3z)}}{\vectornorm{(x,\,y,\,z)^{T}-(x_0,\,y_0,\,z_0)^{T}}}(\text{over ellipsoid}).
\end{align}
It can be seen that each $f_{i}\ (i=1,2,3,4,5)$ stands one class of function.
Function $f_1$, one of Franke functions, was adapted by Renka to the three dimension case \cite{renka1988multivariate}.
$f_1$ is analytic on the sphere.
$f_2$ and $f_3$ were used by Fliege and Maier \cite{fliege1999distribution} to test the quality of their numerical integration scheme, which is based on integration of the polynomial interpolation through their calculated points.
Function $f_2$ , which show in \figref{fig:f2fig}, have only $C^0$ continuity, in particular they are not continuously differentiable at points where any component of $\mathbf{x}$ is zero.
Function $f_3$, which is called ``near-singular function'' \cite{vijayakumar1989new},\ is analytic over $\Stwo$ with a pole just off the surface of the sphere at $\mbfx=(0, 0, 1.01)^{T}$. That is, $f_3((0, 0, 1.01)^{T})=\inf$.
The cosine cap function $f_4$ 
is part of a standard test set for numerical approximations to the shallow water equations in spherical geometry \cite{williamson1992standard}.
$f_4$ is smooth everywhere except at the edge, where two part are joined.
For $f_4$, we set the center $\mathbf{x}_{c}=(x_c,\,y_c,\,z_c)^{T}=(0, 0, 1)^{T}$, radius $R=1/3$ and amplitude $h_0=1$.
Function $f_5$ and $f_6$, used in \cite{sidi2006numericalII} and \cite{atkinson2004quadrature} respectively, are singular functions, which value become infinity at $(x_0,\,y_0,\,z_0)^{T}$.
We make use of above two variable transformations in computation of integrals of these two singular functions.
$f_5$'s singular point is $(0,\,0,\,-1)^{T}$ over $\Stwo$.
The difference between $f_5$ and $f_6$ is that $f_6$ is defined on the ellipsoid
  $$ \mathbb{U} : \left(\frac{x}{1}\right)^2 + \left(\frac{y}{2}\right)^2 + \left(\frac{z}{3}\right)^2=1, $$
and its singular point is $(1/2,\,1,\,3\sqrt{2}/2)^{T}$ over ellipsoid $\mathbb{U}$.
In this case, we assume that a mapping \cite{atkinson2004quadrature}
  \begin{equation}\label{eqt:map_surf}
    \mathcal{M}: \Stwo \xlongrightarrow[\text{onto}]{1-1} \mathbb{U}
  \end{equation}
is given with $\Stwo$.
The integral becames
  \begin{equation*}
    I(f) := \int_{\Stwo} f\left(\mathcal{M}(\mathbf{x})\right) J_{\mathcal{M}}(\mathbf{x}) \domegax,
  \end{equation*}
where $J_{\mathcal{M}}(\mathbf{x})$ is the Jacobian of the mapping $\mathcal{M}$.
With the ellipsoidal surface $\mathbb{U}$ defined as above, we can write
  \begin{equation*}
    \mathcal{M}: (x ,\, y ,\, z)^{T}\in\Stwo \mapsto (\xi ,\, \eta ,\, \zeta)^{T} = (ax ,\, by ,\, cz)^{T}\in\mathbb{U},
  \end{equation*}
and its Jacobian $J_{\mathcal{M}}(\mathbf{x}) = \sqrt{ (bcx)^{2}+(acy)^{2}+(abz)^{2} }$.
This mapping \eqtref{eqt:map_surf} can extend to smooth surface $\mathbb{U}$ which is the boundary of a bounded simply-connected region $\Omega \subset \Real^{3}$ as introduced in \cite{atkinson2004quadrature}.

  \def\funfigRatio{0.47}
  \begin{figure}[htbp]
    \centering
    \subfigure[ $f_{1}$ \label{fig:f1fig}]{
      \includegraphics[width=\funfigRatio\textwidth]{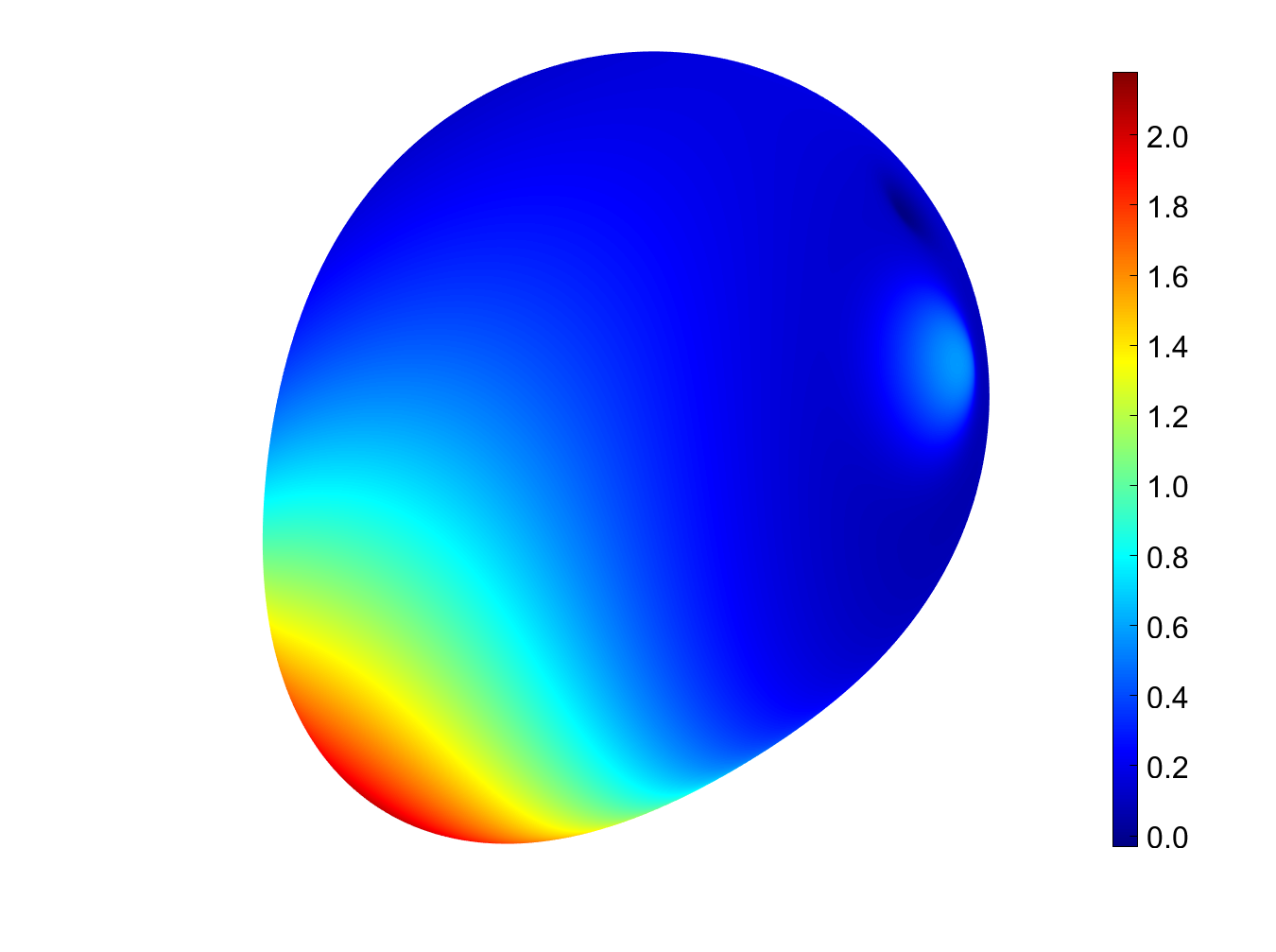}
    }
    \subfigure[ $f_{2}$ \label{fig:f2fig}]{
      \includegraphics[width=\funfigRatio\textwidth]{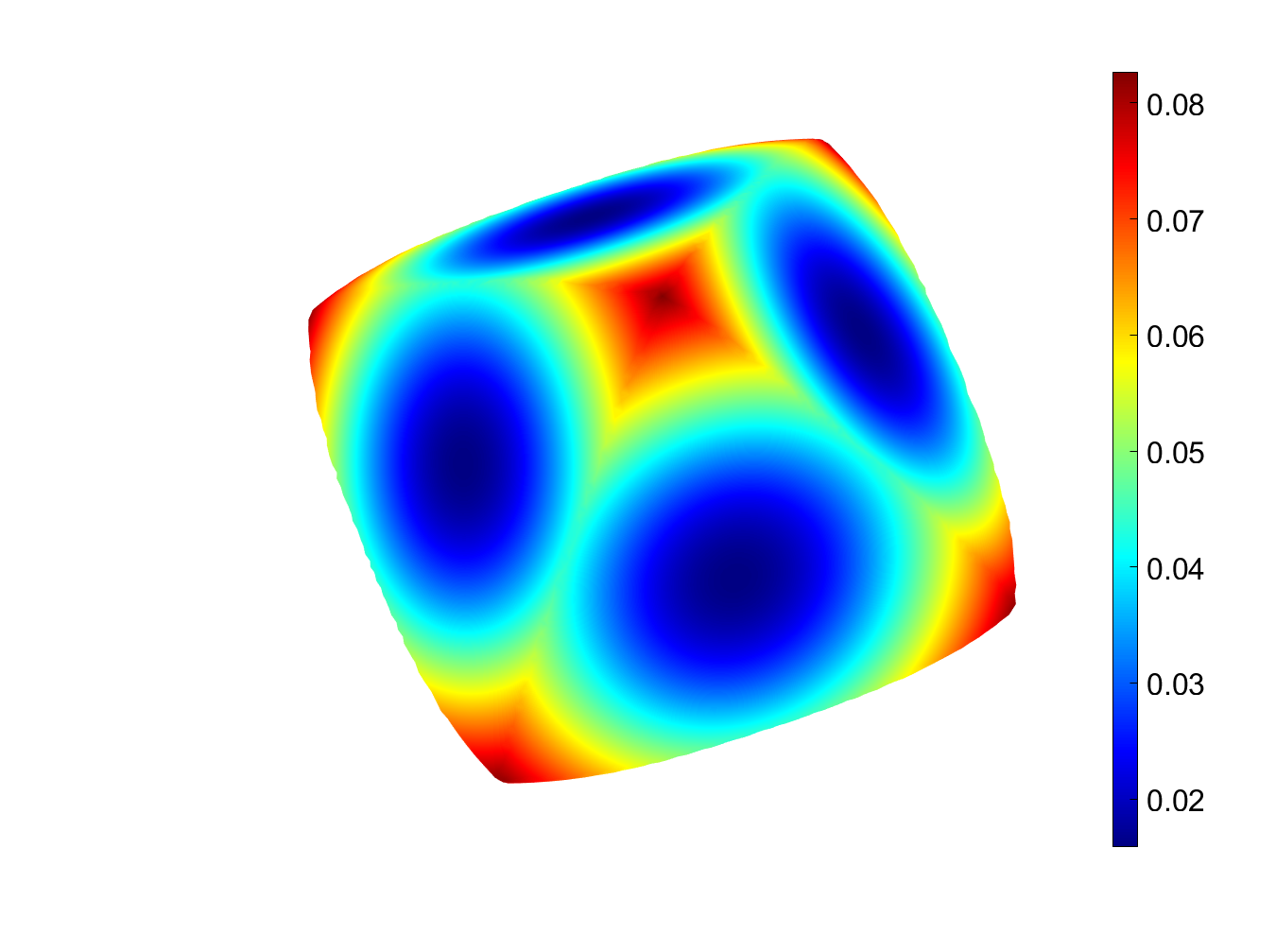}
    }
    \subfigure[ $f_{3}$ \label{fig:f3fig}]{
      \includegraphics[width=\funfigRatio\textwidth]{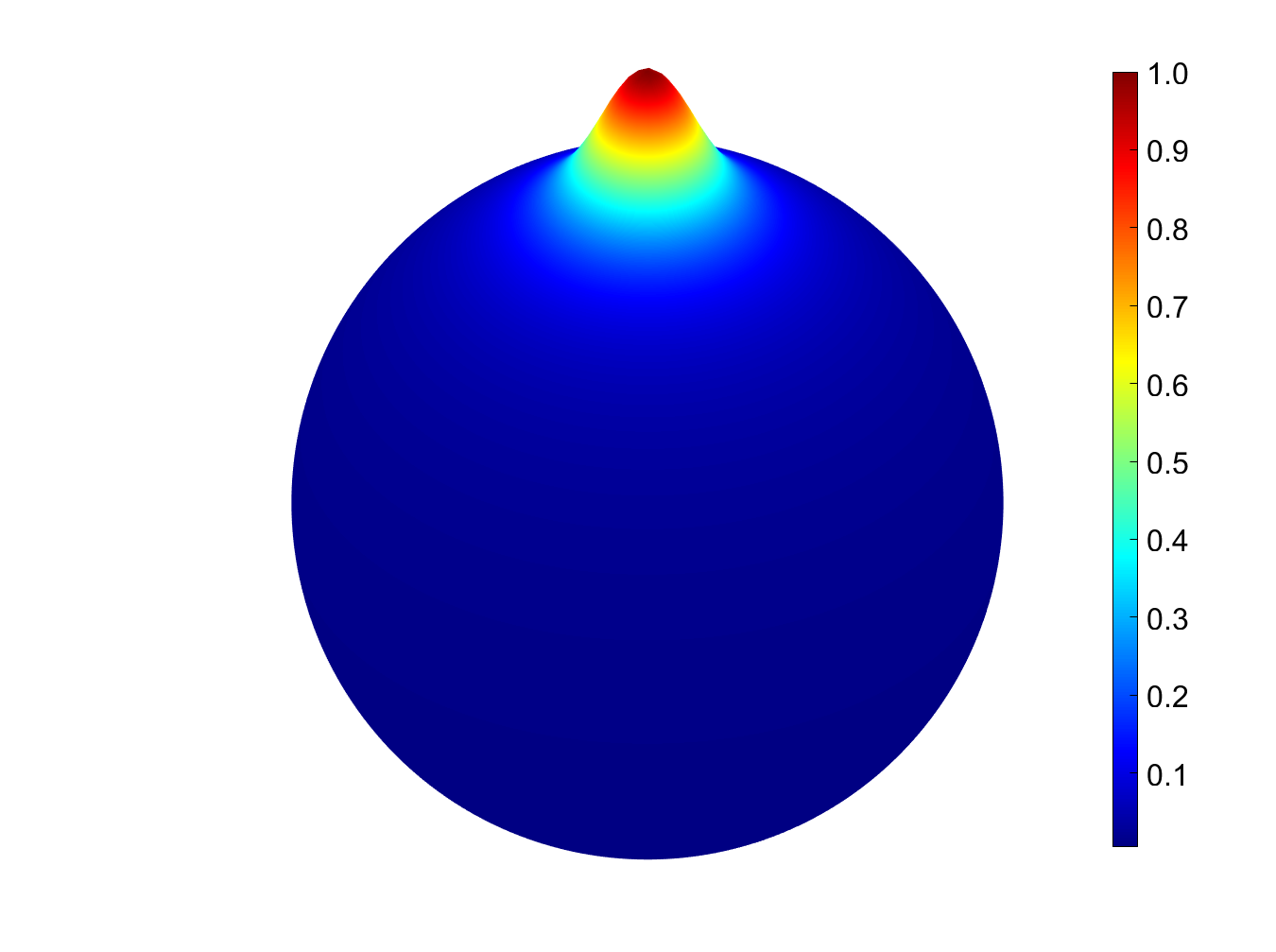}
    }
    \subfigure[ $f_{4}$ \label{fig:f4fig}]{
      \includegraphics[width=\funfigRatio\textwidth]{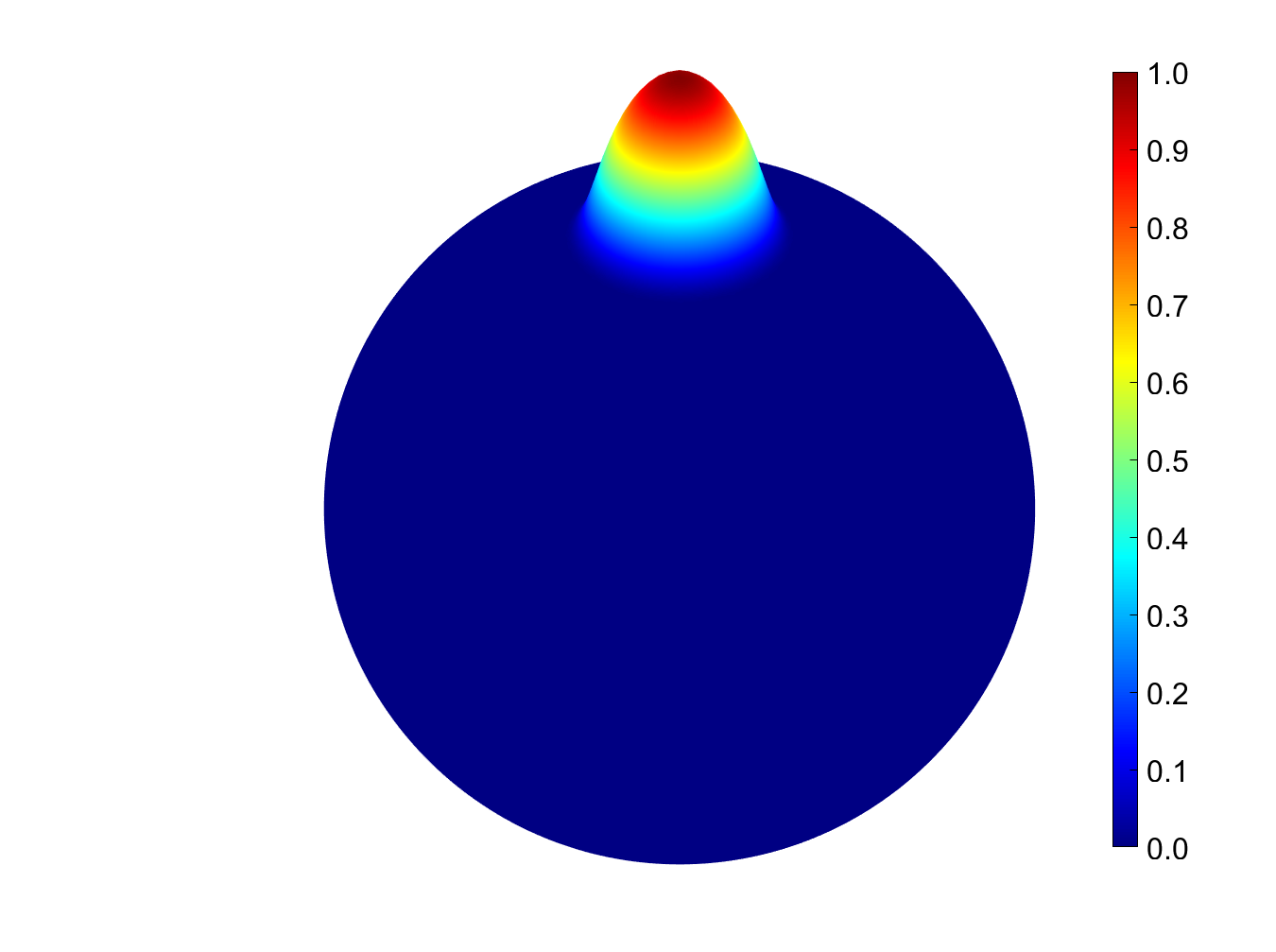}
    }
    \subfigure[ $f_{5}$ \label{fig:f5fig}]{
      \includegraphics[width=\funfigRatio\textwidth]{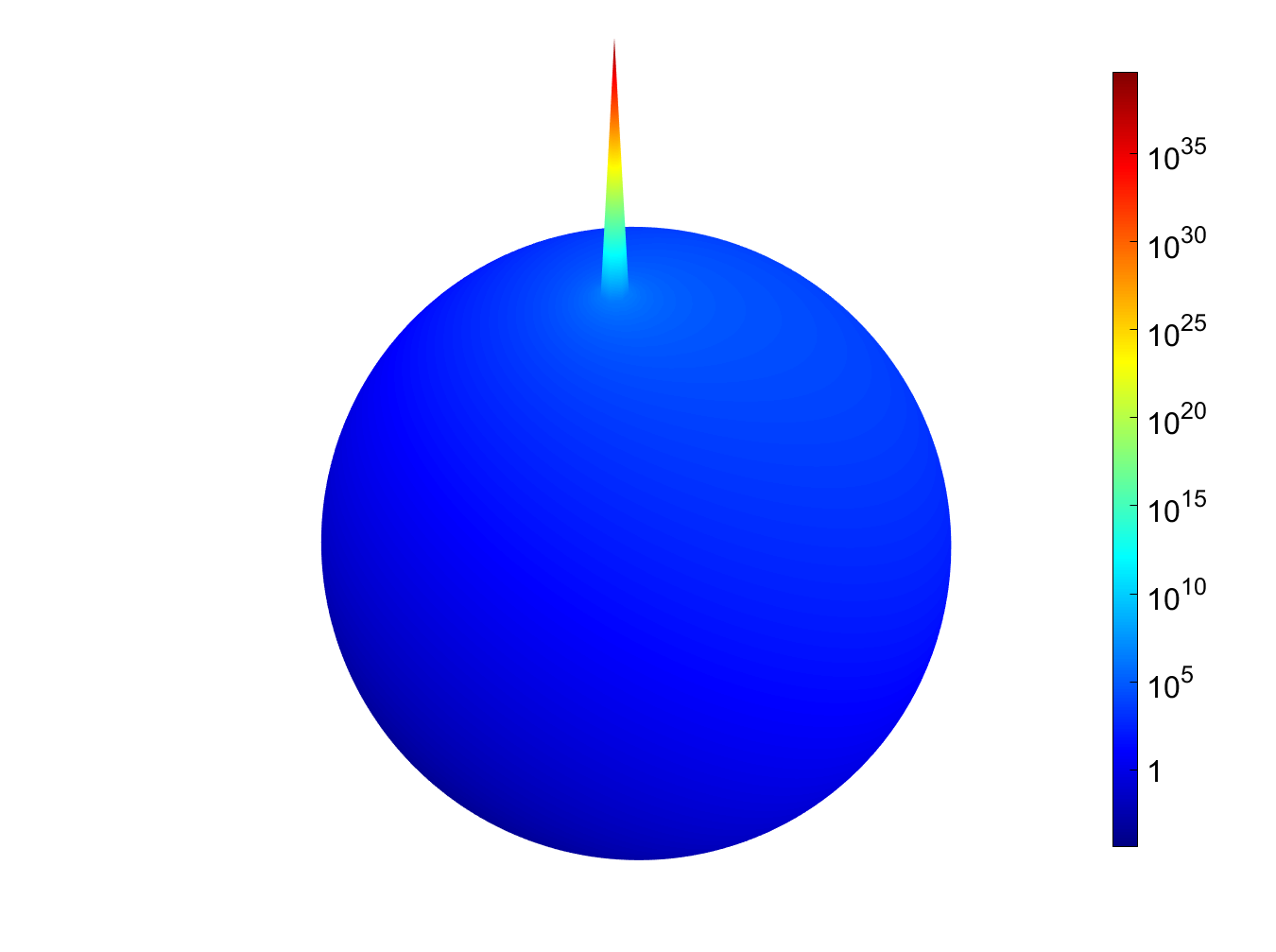}
    }
    \subfigure[ $f_{6}$ \label{fig:f6fig}]{
      \includegraphics[width=\funfigRatio\textwidth]{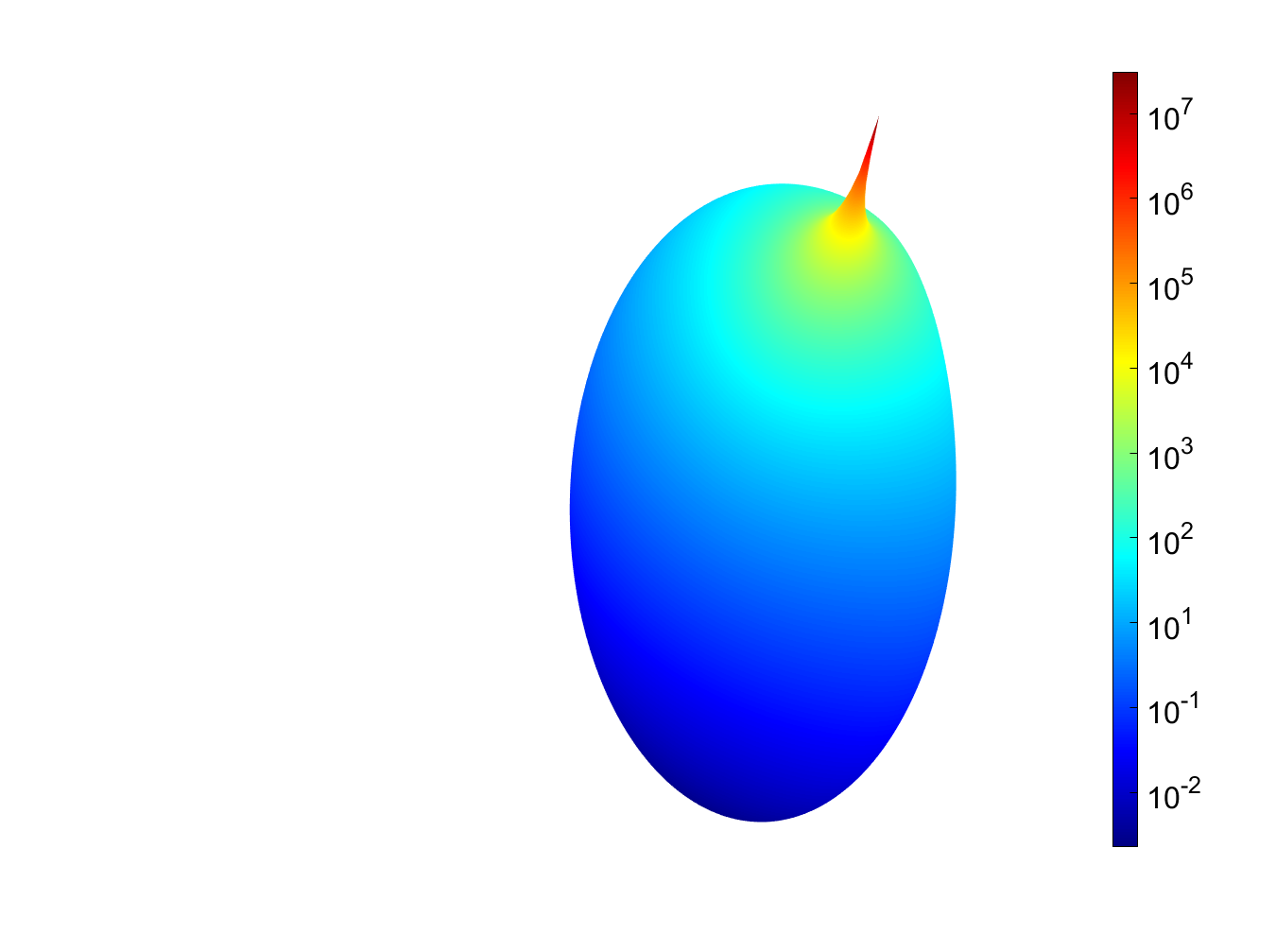}
    }
    \caption{Test functions} \label{fig:funfig}
  \end{figure}

By using Mathematica, exact integration values of all above testing functions over $\Stwo$ are shown in \tabref{tab:intvalue}.  
  \begin{table*}[htbp]
    \centering
      \caption{exact integration value}
      \label{tab:intvalue}
      \begin{tabular}{cl}
        \hline
        function & exact integration values \\
        \hline
        $f_1$    & $6.6961822200736179523\ldots $ \\
        \hline
        $f_2$    & $0.45655373989\ldots$ \\
        \hline
        $f_3$    & $\pi\log{201}/50$ \\
        \hline
        $f_4$    & $0.103351\ldots$ \\
        \hline
        $f_5$    & $40.90220018862976\ldots$ \\
        \hline
        $f_6$    & $371.453416333927\ldots$ \\
        \hline
      \end{tabular}

  \end{table*}

\subsection{Numerical Expertments}

The computational integration error of these six functions by using mentioned quadrature rules are shown in \figref{fig:errAbsAll} .
From \figref{fig:errAbsF1}and \ref{fig:errAbsF2}, it can be seen that WSTD has the best performance in integration when the degree $t$ increases.
The rate of change in integration error of \std{} is sharp as $N$ increases.
For $f_3$ , the Bivariate trapezoidal rule have a bit better than WSTD, but the integration of WSTD also present a competitive descend phenomenon.
In fact, \std{} is a rotationally invariant quadrature rule over $\Stwo$, rather than Bivariate trapezoidal rule depends on latitude and longitude.
For singular functions, we employ Atkinson's transformation and Sidi's transformation.
Then we obtain smoother integrand.
Consequently, the error curve of WSTD performs a rapid descend as $N$ increases, see \figref{fig:errAbsF5_1} \ref{fig:errAbsF5_2}.
It is evident that the error curves of other two quadrature nodes slides slowly even when $N$ passes $10^{4}$. 
For $f_6$ over ellipsoid, the errors slips totally.
But WSTD and Bivariate trapezoidal rule show similar sharp decline phenomenon.
The rate of descent of Equal partition area points is gently as shown by start symbol, see \figref{fig:errAbsF6_2}.

\def\errAbsFigRatio{0.72}
\begin{figure}[htbp]
  \centering
  \subfigure[$f_1$ \label{fig:errAbsF1}]{
    \includegraphics[width=\errAbsFigRatio\textwidth]{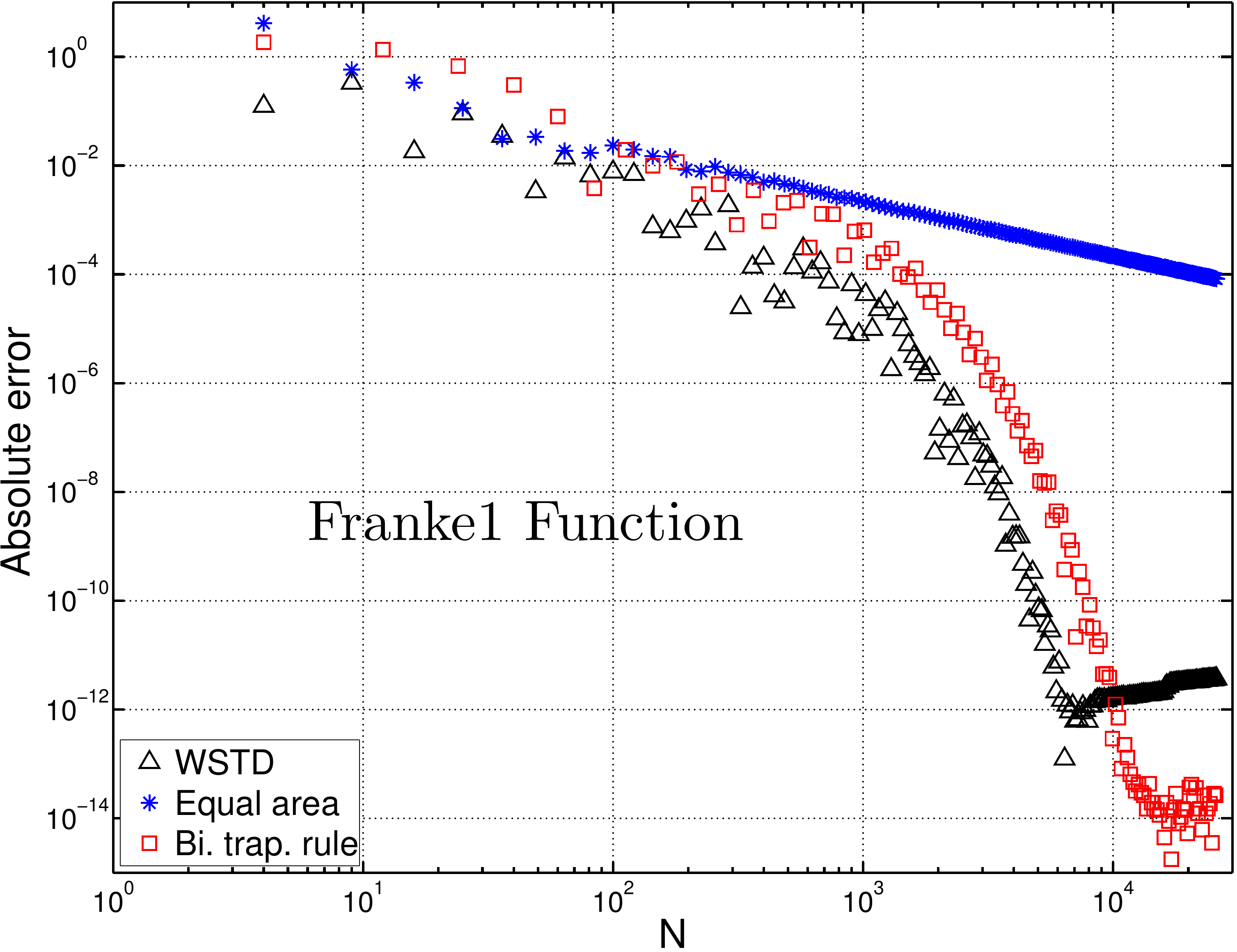}
  }
  \subfigure[$f_2$ \label{fig:errAbsF2}]{
    \includegraphics[width=\errAbsFigRatio\textwidth]{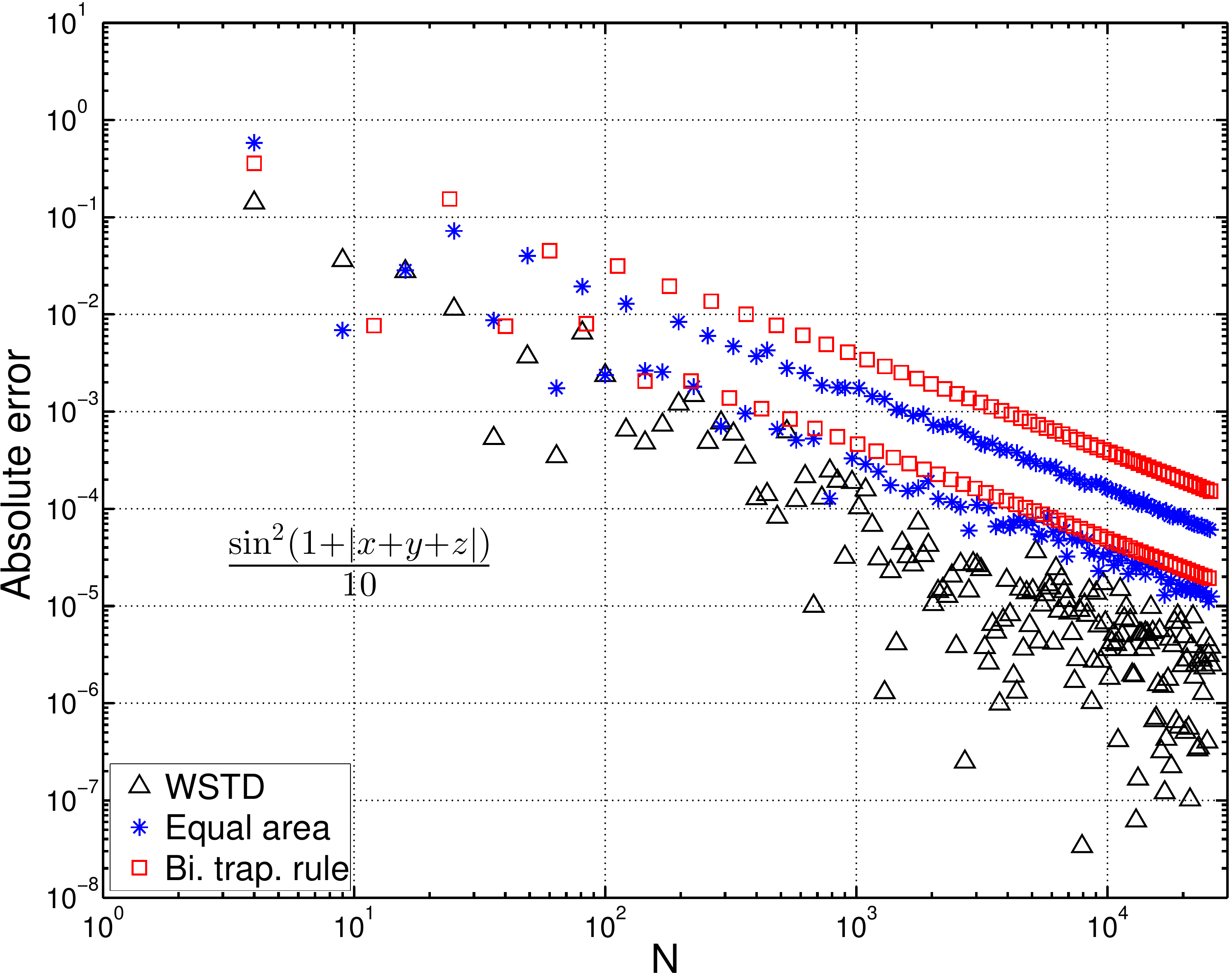}
  }
  \caption{Uniform error for test functions}\label{fig:errAbsAll}
\end{figure}

\begin{figure}[htbp]
  \centering
  \contsubfigure[$f_3$ \label{fig:errAbsF3}]{
    \includegraphics[width=\errAbsFigRatio\textwidth]{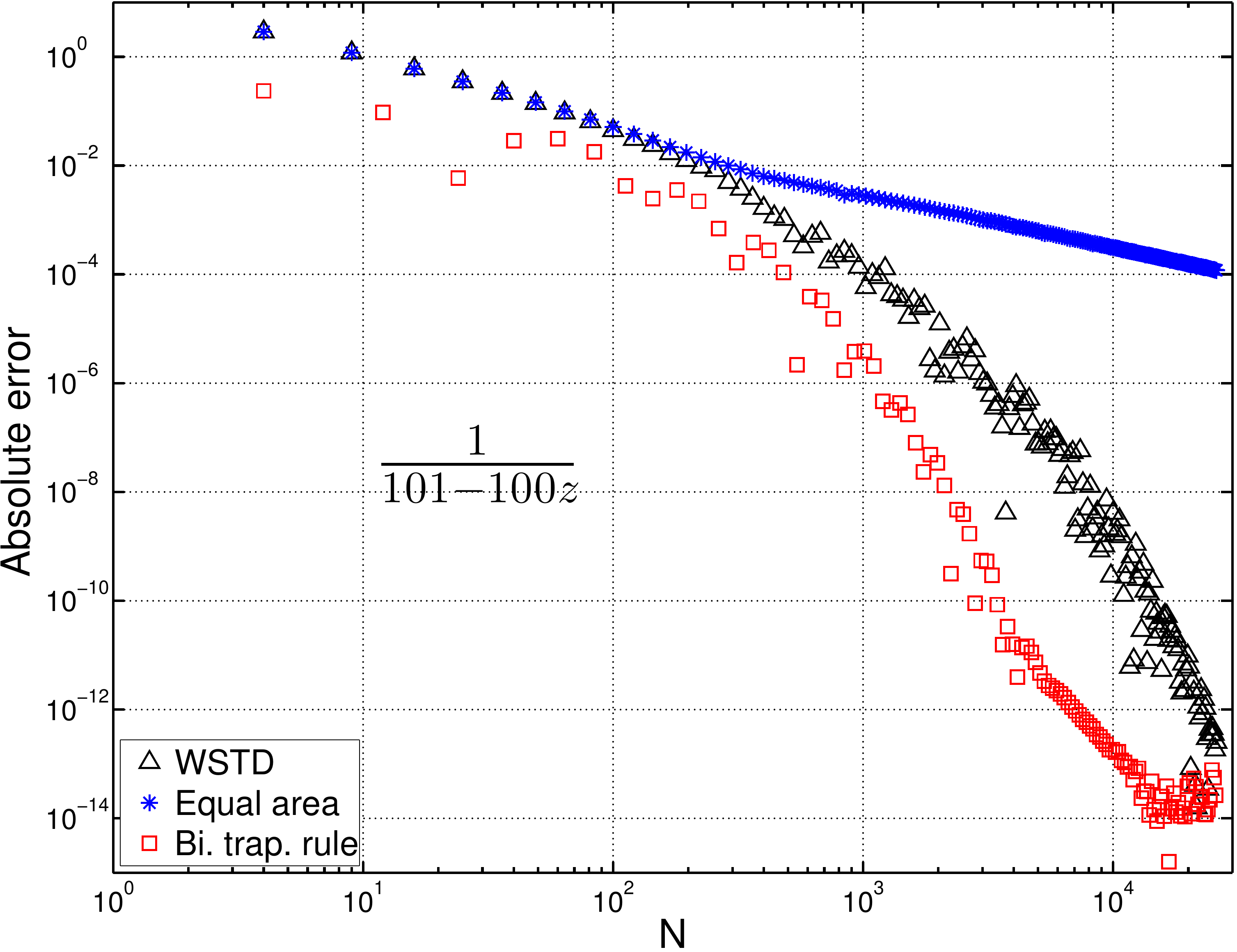}
  }
  \contsubfigure[$f_4$ \label{fig:errAbsF4}]{
    \includegraphics[width=\errAbsFigRatio\textwidth]{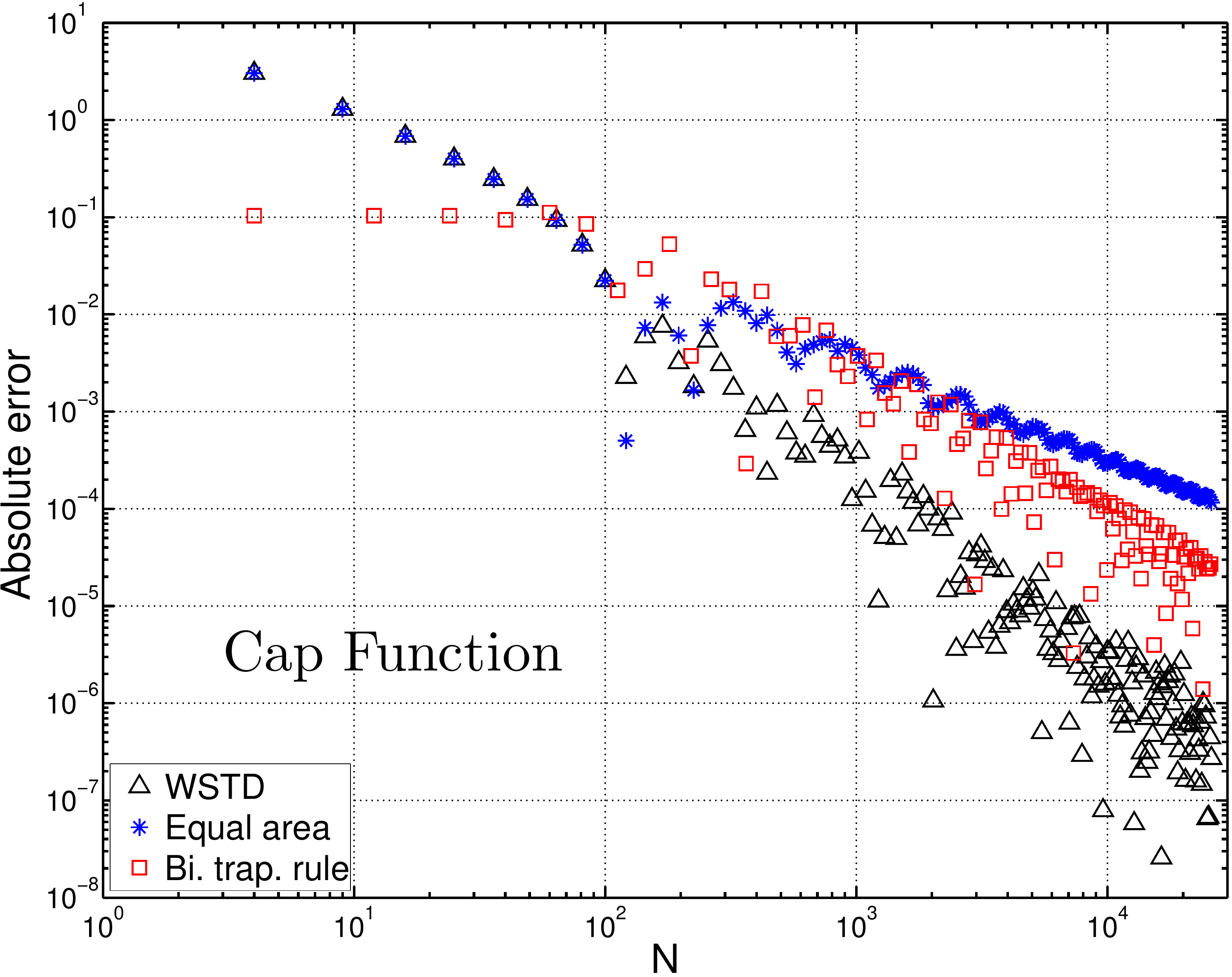}
  }
  \contcaption{Uniform error for test functions}
\end{figure}

\begin{figure}[htbp]
  \centering
  \contsubfigure[$f_5$ with transformation $\Tone$ ( $q=2$ ) \label{fig:errAbsF5_1}]{
    \includegraphics[width=\errAbsFigRatio\textwidth]{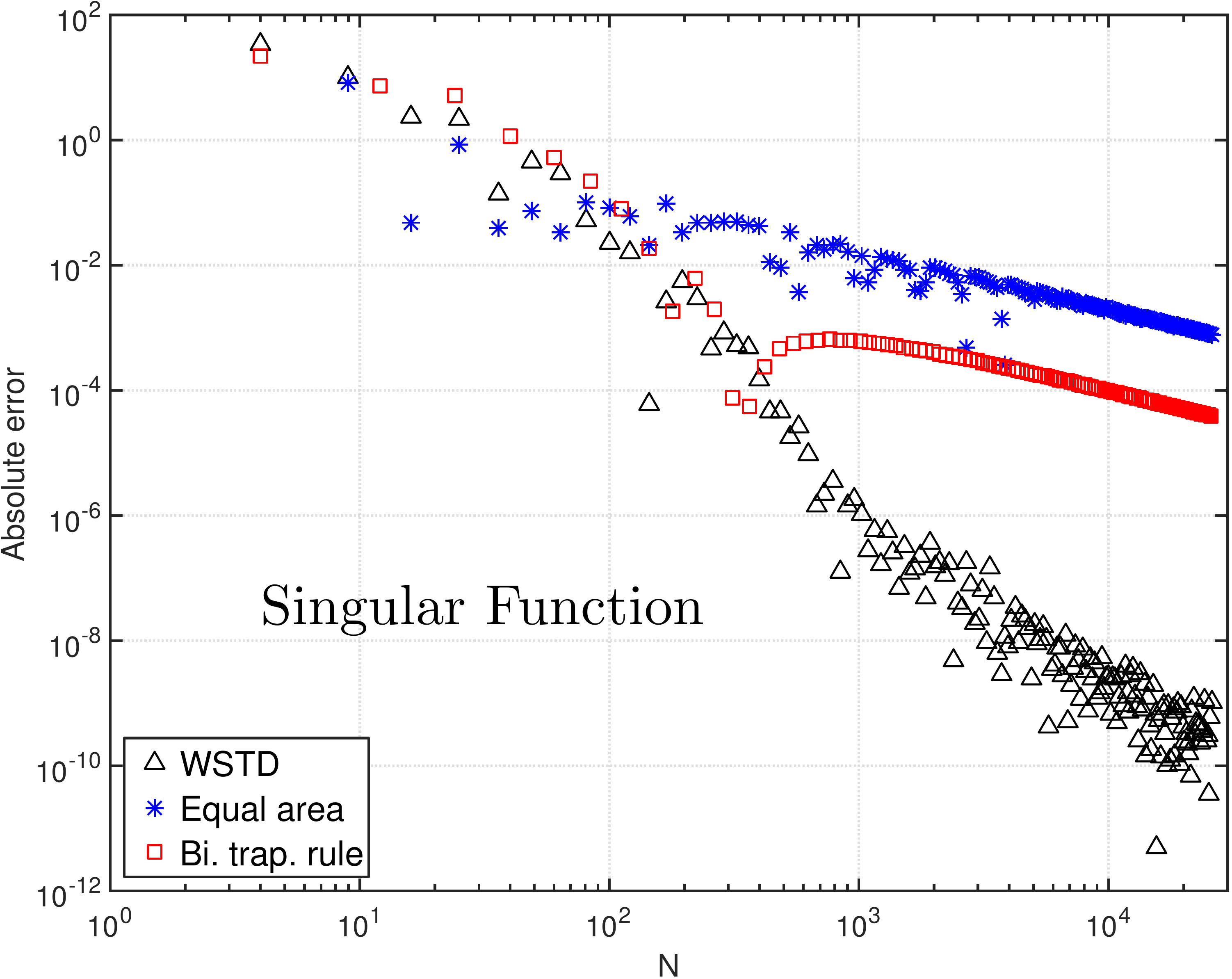}
  }
  \contsubfigure[$f_5$ with transformation $\Ttwo$ ( $m=3$ ) \label{fig:errAbsF5_2}]{
    \includegraphics[width=\errAbsFigRatio\textwidth]{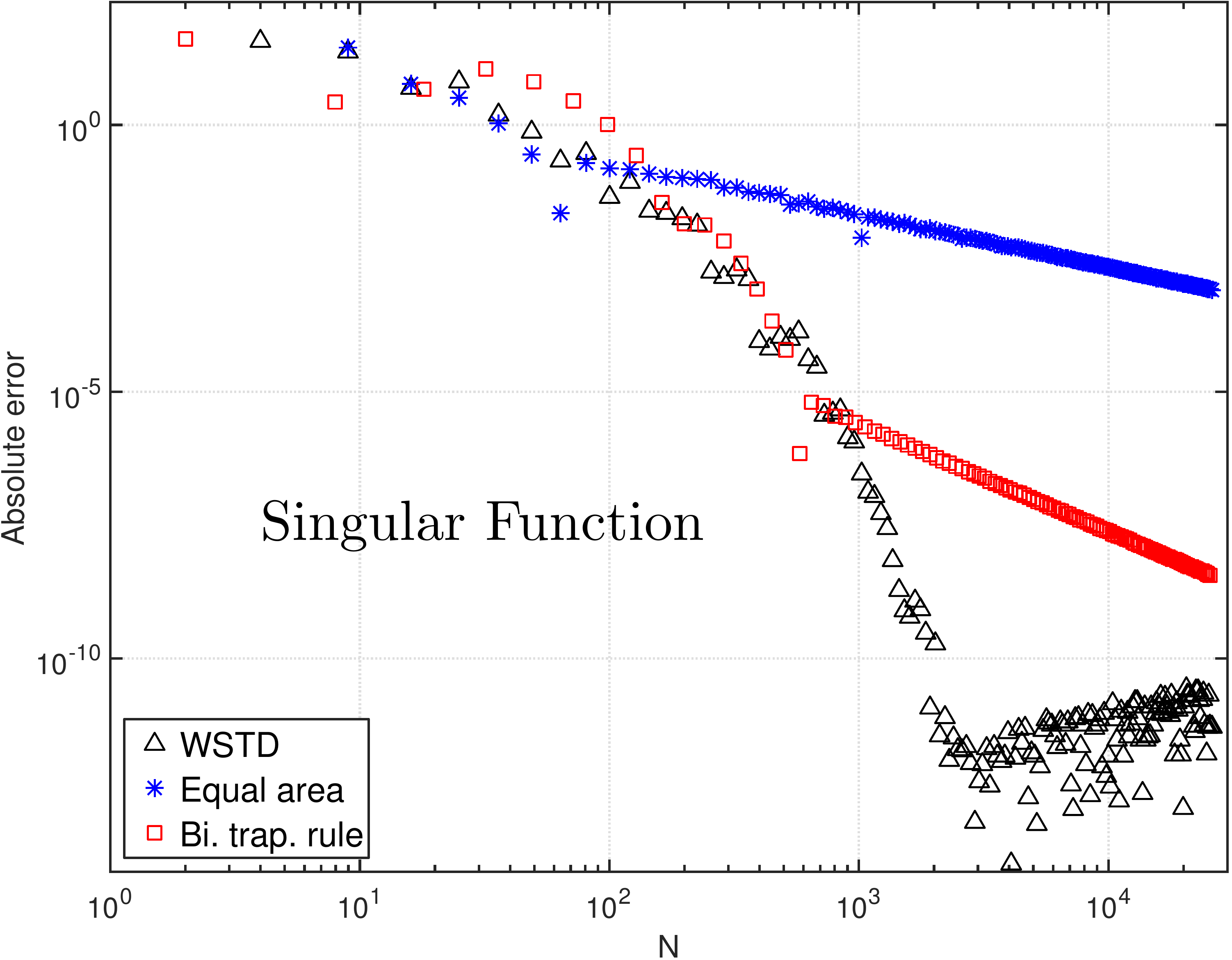}
  }
  \contcaption{Uniform error for test functions}
\end{figure}

\begin{figure}[htbp]
  \centering
  \contsubfigure[$f_6$ with transformation $\Tone$ ( $q=3$ ) \label{fig:errAbsF6_1}]{
    \includegraphics[width=\errAbsFigRatio\textwidth]{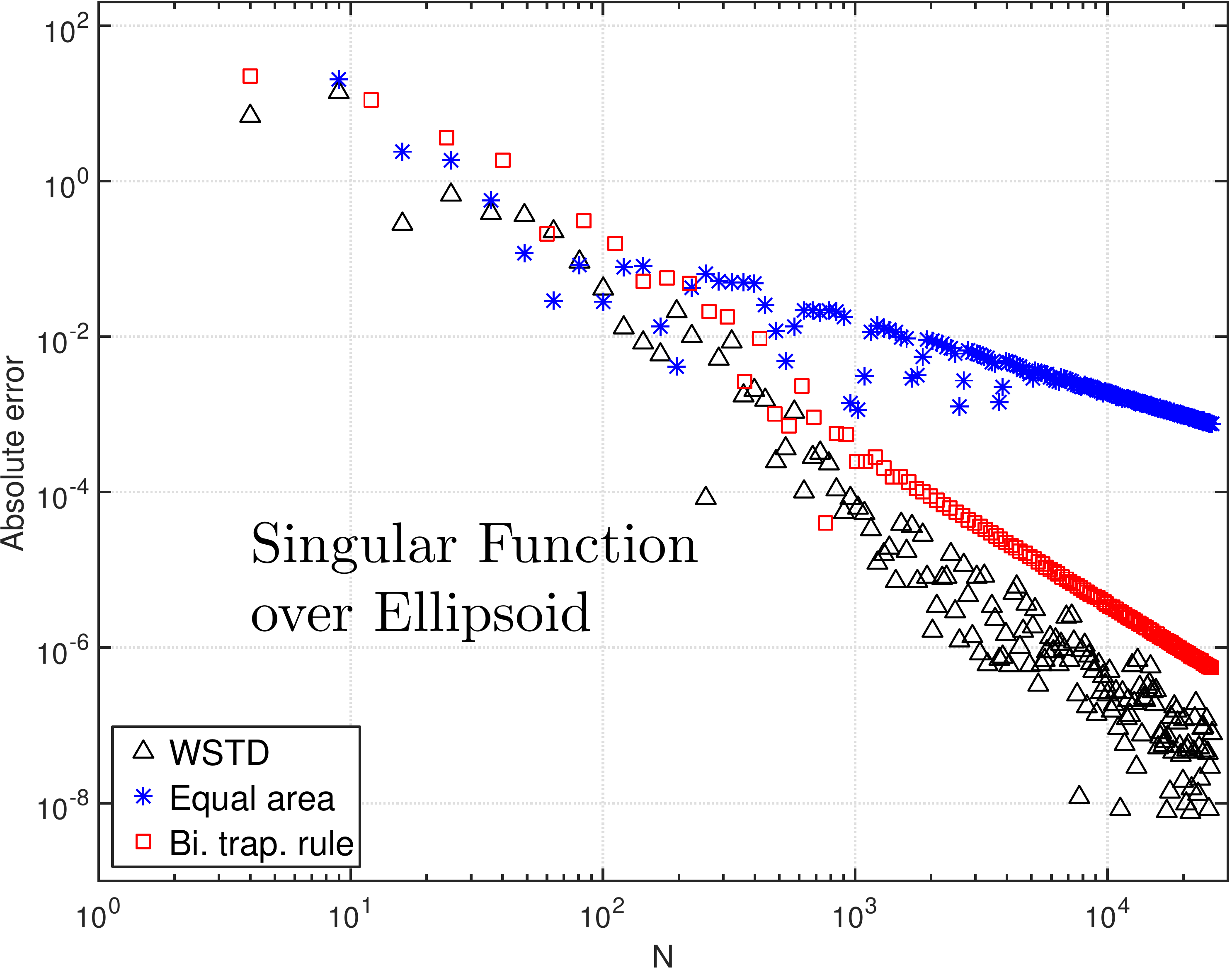}
  }
  \contsubfigure[$f_6$ with transformation $\Ttwo$ ( $m=5$ ) \label{fig:errAbsF6_2}]{
    \includegraphics[width=\errAbsFigRatio\textwidth]{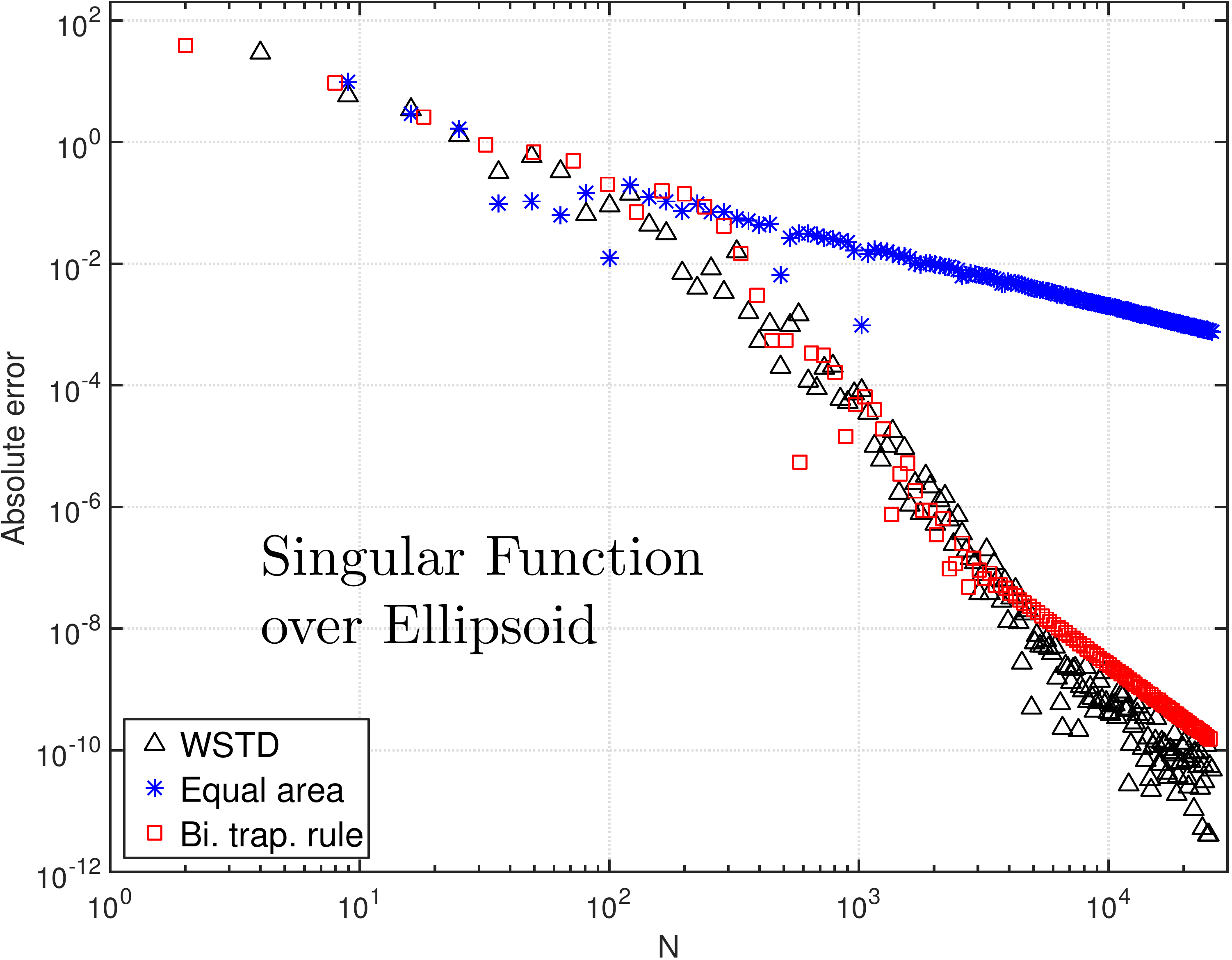}
  }
  \contcaption{Uniform error for test functions}
  \subconcluded
\end{figure}

%% file: chapter/Final-Remark.tex
\section{Final Remark}
The above test results and discussion has been to improve the understanding of properties of quadrature nodes distributions.
We investigated WSTD for approximating the integral of certain functions over the unit sphere, concentrating on the application of WSTD to singular integrands.
By comparison of the computational results of other two quadrature nodes ( Bivariate trapezoidal rule and Equal area points ), WSTD has a remarkable advantage.
All numerical experiments are vivid and encouraging.
Theoretical analysis of these numerical phenomenon is clearly needed in future. 
Further study should be conducted on approximating more complicated integrands over the unit sphere by using WSTD.

%% file: chapter/acknowledgment.tex
\section{Acknowledgment}
The authors thank Professor Kendall E. Atkinson's code in \cite{atkinson2004quadrature}. The support of the National Natural Science Foundation of China (Grant No. 11301222) is gratefully acknowledged. 